\documentclass[11pt]{article}
\usepackage[left=1in,right=1in,top=1in,bottom=1in]{geometry}
\usepackage{times}
\usepackage{expl3}
\usepackage{cite}
\usepackage[table]{xcolor}
\usepackage{multirow}
\usepackage{stackengine} 
\usepackage{hhline}
\usepackage{lipsum}
\usepackage{titlesec}
\usepackage{wrapfig}
\usepackage{epsfig}
\usepackage{graphicx}
\usepackage{amsmath}
\usepackage{amssymb}
\usepackage{epstopdf}
\usepackage{boldline}
\usepackage{calligra}
\usepackage{url}
\usepackage{blindtext}

\newcommand{\define}{\stackrel{\mbox{\tiny def}}{=}}

\newtheorem{theorem}{Theorem}
\newtheorem{corollary}{Corollary}
\newtheorem{lemma}{Lemma}

\usepackage{mathtools}
\usepackage{epstopdf}
\usepackage{balance}
\usepackage{thmtools}
\usepackage{thm-restate}
\usepackage{hyperref}
\usepackage{cleveref}

\usepackage[ruled,vlined]{algorithm2e}
\include{pythonlisting}
\newcommand{\ostar}{\mathbin{\mathpalette\make@circled\star}}

\makeatletter
\newcommand{\removelatexerror}{\let\@latex@error\@gobble}
\makeatother
\makeatletter
\newcommand*{\rom}[1]{\expandafter\@slowromancap\romannumeral #1@}
\makeatother

\ExplSyntaxOn
\newcommand\latinabbrev[1]{
  \peek_meaning:NTF . {
    #1\@}%
  { \peek_catcode:NTF a {
      #1.\@ }%
    {#1.\@}}}
\ExplSyntaxOff

\numberwithin{equation}{section}

\titleclass{\subsubsubsection}{straight}[\subsubsection]

\begin{document}
\vspace{1cm}
\title{T-product Tensors—Part II: Tail Bounds for Sums of Random T-product Tensors}\vspace{1.8cm}
\author{Shih~Yu~Chang\thanks{Shih Yu Chang is with the Department of Applied Data Science,
		San Jose State University, San Jose, CA, U S A. (e-mail: {\tt
			shihyu.chang@sjsu.edu}).} \and Yimin Wei
\thanks{Corresponding author.
  Yimin Wei is with the School of Mathematical Sciences and Shanghai Key Laboratory of Contemporary Applied Mathematics, Fudan University, Shanghai, 200433, PR China. This author is supported in part by the National Natural Science Foundation of China under grants 11771099 and Innovation
  Program of Shanghai Municipal Education Commission.
            (e-mail: {\tt ymwei@fudan.edu.cn, yimin.wei@gmail.com}).
           }}
\maketitle

\begin{abstract}
This paper is the Part II of a serious work about T-product tensors focusing at establishing new probability bounds for sums of random, independent, T-product tensors. These probability bounds characterize large-deviation behavior of the extreme eigenvalue of the sums of random T-product tensors. We apply Lapalace transform method and Lieb's concavity theorem for T-product tensors obtained from our Part I paper, and apply these tools to generalize the classical bounds associated with the names Chernoff, and Bernstein from the scalar to the T-product tensor setting. Tail bounds for the norm of a sum of random rectangular T-product tensors are also derived from corollaries of random Hermitian T-product tensors cases. The proof mechanism is also applied to T-product tensor-valued martingales and T-product tensor-based Azuma, Hoeffding and McDiarmid inequalities are derived. 
\end{abstract}

\begin{keywords}
random T-product tensors, T-product tensor Chernoff bound, T-product tensor Bernstein bound, T-product tensor-valued martingale, T-product tensor Azuma inequality, T-product tensor McDiarmid inequality.
\end{keywords}

\section{Introduction}\label{sec:Introduction} 

\subsection{From Sums of Random Matrices to Sums of Random T-product Tensors}\label{subsec:From Sums of Random Matrices to Sums of Random T-product Tensors} 

In probability theory and theoretical physics, a random matrix is a matrix-valued random variable—that is, a matrix with all entries as random variables. Many crucial physical phenomena can be modeled as random matrix problems. For example, random matrices were introduced by Eugene Wigner to model the nuclei of heavy atoms in nuclear physics~\cite{wigner1993characteristic}. Since then, random matrices have become ubiquitous in science and engineering applications. As this trend accelerates, more and more researchers have to integrate concepts from random matrices into their work. Classical random matrix theory can be difficult to apply, and it is necessary to invent new tools that are easy to use and that apply to a wide range of random matrices~\cite{tropp2019matrix}. Tail bounds for sums of random matrices are among the most popular of these new tools. Tail bounds for sums of random matrices have already found various applications in science and engineering, including: combinatorics~\cite{oliveira2009spectrum}, numerical linear algebra~\cite{martinsson2020randomized}, optimization~\cite{cheung2012linear}, signal processing~\cite{chen2014coherent}, and machine learning~\cite{lopez2014randomized}, etc.

The T-product operation between two three order tensors was introduced by Kilmer and her collaborators in~\cite{kilmer2011factorization, kilmer2013third} to generalize the traditional matrix product. T-product operation has been demonstrated as an important mathematical framework in many fields: multilinear algebra~\cite{li2020continuity, zheng2021t, miao2021t, miao2020generalized}, numerical linear algebra~\cite{zhang2018randomized}, signal processing~\cite{zhang2016exact, semerci2014tensor}, machine learning~\cite{settles2007multiple}, image processing~\cite{khalil2021efficient}, computer vision~\cite{zhang2014novel, martin2013order}, low-rank tensor approximation~\cite{xu2013parallel, zhou2017tensor, qi2021tsingular} etc. However, all these applications assume that systems modelled by T-product tensors are deterministic and such assumption is not true and practical in solving T-product tensors associated issues. In recent years, there are more works begin to study random tensors, see~\cite{chang2020convenient},~\cite{chang2021general},~\cite{chang2021tensor},~\cite{vershynin2020concentration} and references therein. 

In our Part I paper \cite{chang2021T-product}, we  establish following inequalities about T-product tensors: (1) trace function
nondecreasing/convexity; (2) Golden-Thompson inequality for T-product tensors; (3) Jensen’s T-product
inequality; (4) Klein’s T-product inequality. All these inequalities are used to generalize celebrated Lieb’s
concavity theorem from matrices to T-product tensors.

In this work, we will focus on establishing several new tail bounds for sums of random T-product tensors.

\subsection{Tail Bounds Derived in This Paper}\label{subsec:Main Tail Bounds Derived at This Paper}

In this introduction section, we will highlight theorems about tail bounds for sums of random T-product tensors established in this paper. There are two categories of tail bounds discussed here: bounds for eigenvalue and bounds for eigentuples. For bounds related to eigntuples, there is a special condition to be satisfied for the T-product tensor whose eigentuple tail behavior is our interest. 

Let $\mathcal{Y} \in \mathbb{C}^{m \times m \times p}$ be a random T-positive definite (TPD) tensor and we say the tensor $\mathcal{Y}$ satisfies Eq.~\eqref{eq1:lma: Laplace Transform Method Eigentuple Version} if the following inequality relation is valid for the tensor $\mathcal{Y}$:
\begin{eqnarray}\label{eq1:lma: Laplace Transform Method Eigentuple Version}
\frac{1}{p} \lambda_{\max}^{p}(e^{\mathcal{Y}}) + 1 - \frac{1}{p} \leq \mathrm{Tr}(e^{ \mathcal{Y}}), 
\end{eqnarray}
where $t > 0$. If we scale the random TPD tensor $\mathcal{Y}$ as the $\lambda_{\max}(e^{\mathcal{Y}}) = 1$, then Eq.~\eqref{eq1:lma: Laplace Transform Method Eigentuple Version} always holds. 

\subsubsection{Tail Bounds for Sum of Hermitian T-product Tensors with Random Series}

We extend normal-type tail bounds from scalers with  Gaussian and Rademacher random series to T-product tensors with Gaussian and Rademacher random series. The tail bound for the maximum eigenvalue for the sum of Hermitian T-product tensors with Gaussian and Rademacher series is provided by the following Theorem~\ref{thm:TensorGaussianNormalSeries eigenvalue}.

\begin{restatable}[Hermitian T-product Tensor with Gaussian and Rademacher Series Eigenvalue Version]{thm}{TensorGaussianNormalSeriesEigenvalue}\label{thm:TensorGaussianNormalSeries eigenvalue}
Given a finite sequence of fixed T-product tensors $\mathcal{A}_i \in \mathbb{C}^{m \times m \times p}$, and let $\{ \alpha_i \}$ be a finite sequence of independent standard normal variables. We define  
\begin{eqnarray}\label{eq:4_5}
\sigma^2 &\define& \left\Vert \sum\limits_i^{n} \mathcal{A}^2_i \right\Vert,
\end{eqnarray}
then, for all $\theta \geq 0$, we have 
\begin{eqnarray}\label{eq:4_3}
\mathrm{Pr}\left( \lambda_{\max} \left( \sum\limits_{i=1}^{n}\alpha_i \mathcal{A}_i \right) \geq \theta \right) \leq mp e^{-\frac{\theta^2 }{2 \sigma^2}}.
\end{eqnarray}
We use $\left\Vert \mathcal{X} \right\Vert$ for the spectral norm, which is the largest singular value for the T-product tensor $\mathcal{X}$. Then, we have 
\begin{eqnarray}\label{eq:4_4}
\mathrm{Pr}\left( \left\Vert \sum\limits_{i=1}^{n}\alpha_i \mathcal{A}_i \right\Vert \geq \theta \right) \leq 2mp e^{-\frac{\theta^2 }{2 \sigma^2}}.
\end{eqnarray}
This theorem is also valid for a finite sequence of independent Rademacher random variables $\{ \alpha_i \}$.
\end{restatable}

The eigentuple version for T-product tensors with Gaussian and Rademacher random series is provided by the folloiwing Theorem~\ref{thm:TensorGaussianNormalSeries eigentuple}. We use $\left\Vert \mathcal{C} \right\Vert_{\mbox{\tiny{vec}}}$ to represent the spectral norm of eigentuple of the tensor $\mathcal{C}$, which is defined as 
\begin{eqnarray}
\left\Vert \mathcal{C} \right\Vert_{\mbox{\tiny{vec}}} \define \mathbf{d}_{\max}\left( \sqrt{ \mathcal{C}^{\mathrm{H}}\star \mathcal{C}  } \right).
\end{eqnarray}

\begin{restatable}[Hermitian T-product Tensor with Gaussian and Rademacher Series Eigentuple Version]{thm}{TensorGaussianNormalSeriesEigentuple}\label{thm:TensorGaussianNormalSeries eigentuple}
Given a finite sequence of Hermitian T-product tensors $\mathcal{A}_i \in \mathbb{C}^{m \times m \times p}$, and let $\{ \alpha_i \}$ be a finite sequence of independent standard normal variables. We define  
\begin{eqnarray}\label{eq:4_5 eigentuple}
\sigma^2 &\define& \left\Vert \sum\limits_i^{n} \mathcal{A}^2_i \right\Vert,
\end{eqnarray}
then, for all $\mathbf{b} \geq \mathbf{0}$ and $\sum\limits_{i=1}^{n} t \alpha_i \mathcal{A}_i$ satisfying Eq.~\eqref{eq1:lma: Laplace Transform Method Eigentuple Version} for $t >0$, we have 
\begin{eqnarray}\label{eq:4_3 eigentuple}
\mathrm{Pr}\left( \mathbf{d}_{\max} \left( \sum\limits_{i=1}^{n}\alpha_i \mathcal{A}_i \right) \geq \mathbf{b} \right) \leq mp e^{- \frac{b_{\tilde{j}}^2}{2 \sigma^2}},
\end{eqnarray}
where $\tilde{j} \define \arg \min\limits_j \left\{ b_j \right\}$. 
And
\begin{eqnarray}\label{eq:4_4 eigentuple}
\mathrm{Pr}\left( \left\Vert \sum\limits_{i=1}^{n}\alpha_i \mathcal{A}_i \right\Vert_{\mbox{\tiny{vec}}} \geq \mathbf{b} \right) \leq 2 mp e^{- \frac{b_{\tilde{j}}^2}{2 \sigma^2}}.
\end{eqnarray}
This theorem is also valid for a finite sequence of independent Rademacher random variables $\{ \alpha_i \}$.
\end{restatable}

\subsubsection{Chernoff Inequaltities about T-product Tensors}

Next,  we will extend Chernoff bounds of random variables to random T-product tensors.
\begin{restatable}[T-product Tensor Chernoff Bound I]{thm}{TensorChernoffBoundIEigenvalue}\label{thm:TensorChernoffBoundI}
Consider a sequence $\{ \mathcal{X}_i  \in \mathbb{C}^{m \times m \times p} \}$ of independent, random, Hermitian T-product tensors that satisfy
\begin{eqnarray}
\mathcal{X}_i \succeq \mathcal{O} \mbox{~~and~~} \lambda_{\max}(\mathcal{X}_i) \leq 1 
\mbox{~~ almost surely.}
\end{eqnarray}
Define following two quantaties:
\begin{eqnarray}
\overline{\mu}_{\max} \define \lambda_{\max}\left( \frac{1}{n} \sum\limits_{i=1}^{n} \mathbb{E} \mathcal{X}_i \right) \mbox{~~and~~} 
\overline{\mu}_{\min} \define \lambda_{\min}\left( \frac{1}{n} \sum\limits_{i=1}^{n} \mathbb{E} \mathcal{X}_i \right),
\end{eqnarray}
then, we have following two inequalities:
\begin{eqnarray}\label{eq:Chernoff I Upper Bound}
\mathrm{Pr} \left( \lambda_{\max}\left( \frac{1}{n}\sum\limits_{i=1}^{n} \mathcal{X}_i \right) \geq \theta \right) \leq mp e^{- n \mathfrak{D}(\theta || \overline{\mu}_{\max} )},\mbox{~~ for $\overline{\mu}_{\max} \leq \theta \leq 1$;}
\end{eqnarray}
and
\begin{eqnarray}\label{eq:Chernoff I Lower Bound}
\mathrm{Pr} \left( \lambda_{\min}\left( \frac{1}{n}\sum\limits_{i=1}^{n} \mathcal{X}_i \right) \leq \theta \right) \leq mp e^{- n \mathfrak{D}(\theta || \overline{\mu}_{\min} )},\mbox{~~ for $0 \leq \theta \leq \overline{\mu}_{\min}$.}
\end{eqnarray}
\end{restatable}

The other version of T-product tensor Chernoff bound by changing $\overline{\mu}_{\max} (\overline{\mu}_{\min}) $ to  $\mu_{\max} (\mu_{\min}) $ (without average with respect to the number of T-product tensors) is provided by the following Theorem~\ref{thm:TensorChernoffBoundII}

\begin{restatable}[T-product Tensor Chernoff Bound II]{thm}{TensorChernoffBoundIIEigenvalue}\label{thm:TensorChernoffBoundII}
Consider a sequence $\{ \mathcal{X}_i  \in \mathbb{C}^{m \times m \times p } \}$ of independent, random, Hermitian tensors that satisfy
\begin{eqnarray}
\mathcal{X}_i \succeq \mathcal{O} \mbox{~~and~~} \lambda_{\max}(\mathcal{X}_i) \leq T
\mbox{~~ almost surely.}
\end{eqnarray}
Define following two quantaties:
\begin{eqnarray}
\mu_{\max} \define \lambda_{\max}\left( \sum\limits_{i=1}^{n} \mathbb{E} \mathcal{X}_i \right) \mbox{~~and~~} 
\mu_{\min} \define \lambda_{\min}\left( \sum\limits_{i=1}^{n} \mathbb{E} \mathcal{X}_i \right),
\end{eqnarray}
then, we have following two inequalities:
\begin{eqnarray}\label{eq:Chernoff II Upper Bound}
\mathrm{Pr} \left( \lambda_{\max}\left( \sum\limits_{i=1}^{n} \mathcal{X}_i \right) \geq (1+\theta) \mu_{\max} \right) \leq mp \left(\frac{e^{\theta}}{ (1 + \theta)^{1 + \theta}  }\right)^{\mu_{\max}/T} ,\mbox{~~ for $\theta \geq 0$;}
\end{eqnarray}
and
\begin{eqnarray}\label{eq:Chernoff II Lower Bound}
\mathrm{Pr} \left( \lambda_{\min}\left( \sum\limits_{i=1}^{n} \mathcal{X}_i \right) \leq (1 - \theta) \mu_{\min} \right) \leq mp \left(\frac{e^{-\theta}}{ (1 - \theta)^{1 - \theta}  }\right)^{\mu_{\min}/T} ,\mbox{~~ for $\theta \in [0,1]$.}
\end{eqnarray}
\end{restatable}

Below are theorems about Chernoff bounds for the maximum and the minimum eigentuples. Theorem~\ref{thm:TensorChernoffBoundI for Eigentuple} is correspond to Theorem~\ref{thm:TensorChernoffBoundI}, and Theorem~\ref{thm:TensorChernoffBoundII eigentuple} is correspond to Theorem~\ref{thm:TensorChernoffBoundII}.

\begin{restatable}[T-product Tensor Chernoff Bound I for Eigentuple]{thm}{TensorChernoffBoundIEigentuple}\label{thm:TensorChernoffBoundI for Eigentuple}
Consider a sequence $\{ \mathcal{X}_i  \in \mathbb{C}^{m \times m \times p} \}$ of independent, random, Hermitian T-product tensors that satisfy
\begin{eqnarray}
\mathcal{X}_i \succeq \mathcal{O} \mbox{~~and~~} \lambda_{\max}(\mathcal{X}_i) \leq 1 
\mbox{~~ almost surely.}
\end{eqnarray}
Define following two quantaties:
\begin{eqnarray}
\overline{\mu}_{\max} \define \lambda_{\max}\left( \frac{1}{n} \sum\limits_{i=1}^{n} \mathbb{E} \mathcal{X}_i \right) \mbox{~~and~~} 
\overline{\mu}_{\min} \define \lambda_{\min}\left( \frac{1}{n} \sum\limits_{i=1}^{n} \mathbb{E} \mathcal{X}_i \right),
\end{eqnarray}
then, given a real vector $\mathbf{b} \geq \mathbf{0} \in \mathbb{R}^p$ with $\tilde{j} \define \arg\min\limits_j \{ b_j \}$ and $\frac{1}{n} \sum\limits_{i=1}^n t \mathcal{X}_i $ satisfing Eq.~\eqref{eq1:lma: Laplace Transform Method Eigentuple Version}, we have following two inequalities:
\begin{eqnarray}\label{eq:Chernoff I Upper Bound eigentuple}
\mathrm{Pr} \left( \mathbf{d}_{\max}\left( \frac{1}{n}\sum\limits_{i=1}^{n} \mathcal{X}_i \right) \geq \mathbf{b} \right) \leq mp e^{- n \mathfrak{D}( \frac{b_{\tilde{j}}}{n}  || \overline{\mu}_{\max} )},\mbox{~~ for $\overline{\mu}_{\max} \leq  \frac{b_{\tilde{j}}}{n}  \leq 1$;}
\end{eqnarray}
and
\begin{eqnarray}\label{eq:Chernoff I Lower Bound eigentuple}
\mathrm{Pr} \left( \mathbf{d}_{\min}\left( \frac{1}{n}\sum\limits_{i=1}^{n} \mathcal{X}_i \right) \leq \mathbf{b} \right) \leq mp e^{- n \mathfrak{D}( \frac{b_{\tilde{j}}}{n}  || \overline{\mu}_{\min} )},\mbox{~~ for $0 \leq  \frac{b_{\tilde{j}}}{n} \leq \overline{\mu}_{\min}$.}
\end{eqnarray}
\end{restatable}

\begin{restatable}[T-product Tensor Chernoff Bound II for Eigentuple]{thm}{TensorChernoffBoundIIEigentuple}\label{thm:TensorChernoffBoundII eigentuple}
Consider a sequence $\{ \mathcal{X}_i  \in \mathbb{C}^{m \times m \times p } \}$ of independent, random, Hermitian T-product tensors that satisfy
\begin{eqnarray}
\mathcal{X}_i \succeq \mathcal{O} \mbox{~~and~~} \lambda_{\max}(\mathcal{X}_i) \leq T
\mbox{~~ almost surely.}
\end{eqnarray}
Define following two quantaties:
\begin{eqnarray}
\mu_{\max} \define \lambda_{\max}\left( \sum\limits_{i=1}^{n} \mathbb{E} \mathcal{X}_i \right) \mbox{~~and~~} 
\mu_{\min} \define \lambda_{\min}\left( \sum\limits_{i=1}^{n} \mathbb{E} \mathcal{X}_i \right).
\end{eqnarray}
If $\sum\limits_{i=1}^n t \mathcal{X}_i $ satisfies Eq.~\eqref{eq1:lma: Laplace Transform Method Eigentuple Version}, we have following two inequalities:
\begin{eqnarray}\label{eq:Chernoff II Upper Bound eigentuple}
\mathrm{Pr} \left( \mathbf{d}_{\max}\left( \sum\limits_{i=1}^{n} \mathcal{X}_i \right) \geq (1+\theta) \mu_{\max} \mathbf{1} \right) \leq mp \left(\frac{e^{\theta}}{ (1 + \theta)^{1 + \theta}  }\right)^{\mu_{\max}/T} ,\mbox{~~ for $\theta \geq 0$;}
\end{eqnarray}
and
\begin{eqnarray}\label{eq:Chernoff II Lower Bound eigentuple}
\mathrm{Pr} \left( \mathbf{d}_{\min}\left( \sum\limits_{i=1}^{n} \mathcal{X}_i \right) \leq (1 - \theta) \mu_{\min} \mathbf{1} \right) \leq mp \left(\frac{e^{-\theta}}{ (1 - \theta)^{1 - \theta}  }\right)^{\mu_{\min}/T} ,\mbox{~~ for $\theta \in [0,1]$.}
\end{eqnarray}
\end{restatable}

\subsubsection{Bernstein Inequaltities about T-product Tensors}

For random variables, Bernstein inequalities give the upper tail of a sum of
independent, zero-mean random variables that are either bounded or subexponential. In this paper, we will extend Bernstein bounds for a sum of zero-mean random T-product tensors.  The bounded T-product tensor Bernstein bounds will be given by Theorem ~\ref{thm:Bounded Tensor Bernstein}. 

\begin{restatable}[T-product Tensor Bernstein Bounds with Bounded $\lambda_{\max}$]{thm}{BoundedTensorBernstein}\label{thm:Bounded Tensor Bernstein}
Given a finite sequence of independent Hermitian T-product tensors $\{ \mathcal{X}_i  \in \mathbb{C}^{m \times m \times p} \}$ that satisfy
\begin{eqnarray}\label{eq1:thm:Bounded Tensor Bernstein}
\mathbb{E} \mathcal{X}_i = 0 \mbox{~~and~~} \lambda_{\max}(\mathcal{X}_i) \leq T 
\mbox{~~almost surely.} 
\end{eqnarray}
Define the total varaince $\sigma^2$ as: $\sigma^2 \define \left\Vert \sum\limits_i^n \mathbb{E} \left( \mathcal{X}^2_i \right) \right\Vert$.
Then, we have following inequalities:
\begin{eqnarray}\label{eq2:thm:Bounded Tensor Bernstein}
\mathrm{Pr} \left( \lambda_{\max}\left( \sum\limits_{i=1}^{n} \mathcal{X}_i \right)\geq \theta \right) \leq mp \exp \left( \frac{-\theta^2/2}{\sigma^2 + T\theta/3}\right);
\end{eqnarray}
and
\begin{eqnarray}\label{eq3:thm:Bounded Tensor Bernstein}
\mathrm{Pr} \left( \lambda_{\max}\left( \sum\limits_{i=1}^{n} \mathcal{X}_i \right)\geq \theta \right) \leq mp \exp \left( \frac{-3 \theta^2}{ 8 \sigma^2}\right)~~\mbox{for $\theta \leq \sigma^2/T$};
\end{eqnarray}
and
\begin{eqnarray}\label{eq4:thm:Bounded Tensor Bernstein}
\mathrm{Pr} \left( \lambda_{\max}\left( \sum\limits_{i=1}^{n} \mathcal{X}_i \right)\geq \theta \right) \leq mp \exp \left( \frac{-3 \theta}{ 8 T } \right)~~\mbox{for $\theta \geq \sigma^2/T$}.
\end{eqnarray}
\end{restatable}

Below is the subexponential T-product tensor Bernstein bounds. Different from Theorem~\ref{thm:Bounded Tensor Bernstein}, we relax the bounded constraint for the 
maximum eigenvalue for T-product tensors $\mathcal{X}_i$ to $\mathbb{E} (\mathcal{X}^p_i) \preceq \frac{p! T^{p-2} }{2} \mathcal{A}_i^2$, where $p = 2,3,4,\cdots$. 

\begin{restatable}[Subexponential T-product Tensor Bernstein Bounds]{thm}{SubexponentialTensorBernstein}\label{thm:Subexponential Tensor Bernstein}
Given a finite sequence of independent Hermitian T-product tensors $\{ \mathcal{X}_i  \in \mathbb{C}^{m \times m \times p} \}$ that satisfy
\begin{eqnarray}\label{eq1:thm:Subexponential Tensor Bernstein}
\mathbb{E} \mathcal{X}_i = 0 \mbox{~~and~~} \mathbb{E} (\mathcal{X}^p_i) \preceq \frac{p! T^{p-2} }{2} \mathcal{A}_i^2,
\end{eqnarray}
where $p = 2,3,4,\cdots$. 

Define the total varaince $\sigma^2$ as: $\sigma^2 \define \left\Vert \sum\limits_i^n \mathcal{A}_i^2 \right\Vert$.
Then, we have following inequalities:
\begin{eqnarray}\label{eq2:thm:Subexponential Tensor Bernstein}
\mathrm{Pr} \left( \lambda_{\max}\left( \sum\limits_{i=1}^{n} \mathcal{X}_i \right)\geq \theta \right) \leq mp \exp \left( \frac{-\theta^2/2}{\sigma^2 + T\theta}\right);
\end{eqnarray}
and
\begin{eqnarray}\label{eq3:thm:Subexponential Tensor Bernstein}
\mathrm{Pr} \left( \lambda_{\max}\left( \sum\limits_{i=1}^{n} \mathcal{X}_i \right)\geq \theta \right) \leq mp \exp \left( \frac{-\theta^2}{ 4 \sigma^2}\right)~~\mbox{for $\theta \leq \sigma^2/T$};
\end{eqnarray}
and
\begin{eqnarray}\label{eq4:thm:Subexponential Tensor Bernstein}
\mathrm{Pr} \left( \lambda_{\max}\left( \sum\limits_{i=1}^{n} \mathcal{X}_i \right)\geq \theta \right) \leq mp \exp \left( \frac{- \theta}{ 4 T } \right)~~\mbox{for $\theta \geq \sigma^2/T$}.
\end{eqnarray}
\end{restatable}

Below are theorems about T-product tensor Bernstein bounds for the maximum and the minimum eigentuples. Theorem~\ref{thm:Bounded Tensor Bernstein eigentuple} is correspond to Theorem~\ref{thm:Bounded Tensor Bernstein}, and Theorem~\ref{thm:Subexponential Tensor Bernstein eigentuple} is correspond to Theorem~\ref{thm:Subexponential Tensor Bernstein}.

\begin{restatable}[T-product Tensor Bernstein Bounds with Bounded $\lambda_{\max}$ for Eigentuple]{thm}{BoundedTensorBernsteinEigentuple}\label{thm:Bounded Tensor Bernstein eigentuple}
Given a finite sequence of independent Hermitian T-product tensors $\{ \mathcal{X}_i  \in \mathbb{C}^{m \times m \times p} \}$ that satisfy
\begin{eqnarray}\label{eq1:thm:Bounded Tensor Bernstein eigentuple}
\mathbb{E} \mathcal{X}_i = 0 \mbox{~~and~~} \lambda_{\max}(\mathcal{X}_i) \leq T 
\mbox{~~almost surely.} 
\end{eqnarray}
Define the total varaince $\sigma^2$ as: $\sigma^2 \define \left\Vert \sum\limits_i^n \mathbb{E} \left( \mathcal{X}^2_i \right) \right\Vert$.
Then, given a positive real vector $\mathbf{b} \geq \mathbf{0} \in \mathbb{R}^p$ with $\tilde{j} \define \arg\min\limits_j \{ b_j \}$ and $\sum\limits_{i=1}^n t \mathcal{X}_i $ satisfing Eq.~\eqref{eq1:lma: Laplace Transform Method Eigentuple Version} for any $t >0$, we have following inequalities:
\begin{eqnarray}\label{eq2:thm:Bounded Tensor Bernstein eigentuple}
\mathrm{Pr} \left( \mathbf{d}_{\max}\left( \sum\limits_{i=1}^{n} \mathcal{X}_i \right)\geq \mathbf{b} \right) \leq mp \exp \left( \frac{-b_{\tilde{j}}^2/2}{\sigma^2 + T\theta/3}\right);
\end{eqnarray}
and
\begin{eqnarray}\label{eq3:thm:Bounded Tensor Bernstein eigentuple}
\mathrm{Pr} \left( \mathbf{d}_{\max}\left( \sum\limits_{i=1}^{n} \mathcal{X}_i \right)\geq \mathbf{b} \right) \leq mp \exp \left( \frac{-3 b_{\tilde{j}}^2}{ 8 \sigma^2}\right)~~\mbox{for $b_{\tilde{j}} \leq \sigma^2/T$};
\end{eqnarray}
and
\begin{eqnarray}\label{eq4:thm:Bounded Tensor Bernstein eigentuple}
\mathrm{Pr} \left( \mathbf{d}_{\max}\left( \sum\limits_{i=1}^{n} \mathcal{X}_i \right)\geq \mathbf{b} \right) \leq mp \exp \left( \frac{-3 b_{\tilde{j}}}{ 8 T } \right)~~\mbox{for $b_{\tilde{j}} \geq \sigma^2/T$}.
\end{eqnarray}
\end{restatable}

\begin{restatable}[Subexponential T-product Tensor Bernstein Bounds for Eigentuple]{thm}{SubexponentialTensorBernsteinEigentuple}\label{thm:Subexponential Tensor Bernstein eigentuple}
Given a finite sequence of independent Hermitian T-product tensors $\{ \mathcal{X}_i  \in \mathbb{C}^{m \times m \times p} \}$ that satisfy
\begin{eqnarray}\label{eq1:thm:Subexponential Tensor Bernstein eigentuple}
\mathbb{E} \mathcal{X}_i = 0 \mbox{~~and~~} \mathbb{E} (\mathcal{X}^p_i) \preceq \frac{p! T^{p-2} }{2} \mathcal{A}_i^2,
\end{eqnarray}
where $p = 2,3,4,\cdots$. 

Define the total varaince $\sigma^2$ as: $\sigma^2 \define \left\Vert \sum\limits_i^n \mathcal{A}_i^2 \right\Vert$.
Then, given a positive real vector $\mathbf{b} \in \mathbb{R}^p$ with $\tilde{j} \define \arg\min\limits_j \{ b_j \}$ and $\sum\limits_{i=1}^n t \mathcal{X}_i $ satisfing Eq.~\eqref{eq1:lma: Laplace Transform Method Eigentuple Version} for any $t >0$, we have following inequalities:
\begin{eqnarray}\label{eq2:thm:Subexponential Tensor Bernstein eigentuple}
\mathrm{Pr} \left( \mathbf{d}_{\max}\left( \sum\limits_{i=1}^{n} \mathcal{X}_i \right)\geq \mathbf{b} \right) \leq mp \exp \left( \frac{-b_{\tilde{j}}^2/2}{\sigma^2 + Tb_{\tilde{j}}}\right);
\end{eqnarray}
and
\begin{eqnarray}\label{eq3:thm:Subexponential Tensor Bernstein eigentuple}
\mathrm{Pr} \left( \mathbf{d}_{\max}\left( \sum\limits_{i=1}^{n} \mathcal{X}_i \right)\geq \mathbf{b} \right) \leq mp \exp \left( \frac{-b_{\tilde{j}}^2}{ 4 \sigma^2}\right)~~\mbox{for $b_{\tilde{j}} \leq \sigma^2/T$};
\end{eqnarray}
and
\begin{eqnarray}\label{eq4:thm:Subexponential Tensor Bernstein eigentuple}
\mathrm{Pr} \left( \mathbf{d}_{\max}\left( \sum\limits_{i=1}^{n} \mathcal{X}_i \right)\geq \mathbf{b} \right) \leq mp \exp \left( \frac{- b_{\tilde{j}}}{ 4 T } \right)~~\mbox{for $b_{\tilde{j}}\geq \sigma^2/T$}.
\end{eqnarray}
\end{restatable}

\subsubsection{Inequaltities about T-product Tensor Martingales}

T-product tensor Azuma and McDiarmid inequalities will be provided for the maximum eigenvalue and the maximum eigentuple versions. 

\begin{restatable}[T-product Tensor Azuma Inequality for Eigenvalue]{thm}{TensorAzuma}\label{thm:Tensor Azuma}
Given a finite adapted sequence of Hermitian tensors $\{ \mathcal{X}_i  \in \mathbb{C}^{m \times m \times p } \}$ and a fixed sequence of Hermitian T-product tensors $\{ \mathcal{A}_i \}$ that satisfy
\begin{eqnarray}\label{eq1:thm:Tensor Azuma}
\mathbb{E}_{i-1} \mathcal{X}_i = 0 \mbox{~~and~~} \mathcal{X}^2_i \preceq  \mathcal{A}^2_i ~~ \mbox{almost surely}, 
\end{eqnarray}
where $i = 1,2,3,\cdots$. 

Define the total varaince $\sigma^2$ as: $\sigma^2 \define \left\Vert \sum\limits_i^n \mathcal{A}_i^2 \right\Vert$.
Then, we have following inequalities:
\begin{eqnarray}\label{eq2:thm:Tensor Azuma}
\mathrm{Pr} \left( \lambda_{\max}\left( \sum\limits_{i=1}^{n} \mathcal{X}_i \right)\geq \theta \right) \leq mp e^{-\frac{\theta^2}{8 \sigma^2}}.
\end{eqnarray}
\end{restatable}

\begin{restatable}[T-product Tensor McDiarmid Inequality]{thm}{TensorMcDiarmid}
\label{thm:Tensor McDiarmid}
Given a set of $n$ independent random variables, i.e. $\{X_i: i = 1,2,\cdots n\}$, and let 
$F$ be a Hermitian T-product tensor-valued function that maps these $n$ random variables to a Hermitian T-product tensor of dimension within $\mathbb{C}^{m \times m \times p}$. Consider a sequence of Hermitian tensors $\{ \mathcal{A}_i \}$ that satisfy
\begin{eqnarray}\label{eq1:thm:Tensor McDiarmid}
\left( F(x_1,\cdots,x_i,\cdots,x_n)  -  F(x_1,\cdots,x'_i,\cdots,x_n) \right)^2 \preceq \mathcal{A}^2_i,
\end{eqnarray}
where $x_i, x'_i \in X_i$ and $1 \leq i \leq n$. Define the total variance $\sigma^2$ as: $\sigma^2 \define \left\Vert \sum\limits_i^n \mathcal{A}_i^2 \right\Vert $.
Then, we have following inequality:
\begin{eqnarray}\label{eq2:thm:Tensor McDiarmid}
\mathrm{Pr} \left( \lambda_{\max}\left( F(x_1,\cdots,x_n) - \mathbb{E}F(x_1,\cdots,x_n )\right)\geq \theta \right) \leq mp  e^{-\frac{\theta^2}{8 \sigma^2}}.
\end{eqnarray}
\end{restatable}

Following two theorems are eigentuple version for T-product tensor Azuma and McDiarmid inequalities. 

\begin{restatable}[T-product Tensor Azuma Inequality for Eigentuple]{thm}{TensorAzumaEigentuple}\label{thm:Tensor Azuma eigentuple}
Given a finite adapted sequence of Hermitian tensors $\{ \mathcal{X}_i  \in \mathbb{C}^{m \times m \times p } \}$ and a fixed sequence of Hermitian T-product tensors $\{ \mathcal{A}_i \}$ that satisfy
\begin{eqnarray}\label{eq1:thm:Tensor Azuma}
\mathbb{E}_{i-1} \mathcal{X}_i = 0 \mbox{~~and~~} \mathcal{X}^2_i \preceq  \mathcal{A}^2_i ~~ \mbox{almost surely}, 
\end{eqnarray}
where $i = 1,2,3,\cdots$. 

Define the total varaince $\sigma^2$ as: $\sigma^2 \define \left\Vert \sum\limits_i^n \mathcal{A}_i^2 \right\Vert$.
Then, given a positive real vector $\mathbf{b} \in \mathbb{R}^p$ with $\tilde{j} \define \arg\min\limits_j \{ b_j \}$ and $\sum\limits_{i=1}^n t \mathcal{X}_i $ satisfing Eq.~\eqref{eq1:lma: Laplace Transform Method Eigentuple Version} for any $t >0$, we have following inequalities:
\begin{eqnarray}\label{eq2:thm:Tensor Azuma}
\mathrm{Pr} \left( \mathbf{d}_{\max}\left( \sum\limits_{i=1}^{n} \mathcal{X}_i \right)\geq \mathbf{b} \right) \leq mp e^{-\frac{b_{\tilde{j}}^2}{8 \sigma^2}}.
\end{eqnarray}
\end{restatable}

\begin{restatable}[T-product Tensor McDiarmid Inequality for Eigentuple]{thm}{TensorMcDiarmidEigentuple}\label{thm:Tensor McDiarmid eigentuple}
Given a set of $n$ independent random variables, i.e. $\{X_i: i = 1,2,\cdots n\}$, and let 
$F$ be a Hermitian T-product tensor-valued function that maps these $n$ random variables to a Hermitian T-product tensor of dimension within $\mathbb{C}^{m \times m \times p}$. Consider a sequence of Hermitian tensors $\{ \mathcal{A}_i \}$ that satisfy
\begin{eqnarray}\label{eq1:thm:Tensor McDiarmid eigentuple}
\left( F(x_1,\cdots,x_i,\cdots,x_n)  -  F(x_1,\cdots,x'_i,\cdots,x_n) \right)^2 \preceq \mathcal{A}^2_i,
\end{eqnarray}
where $x_i, x'_i \in X_i$ and $1 \leq i \leq n$. Define the total variance $\sigma^2$ as: $\sigma^2 \define \left\Vert \sum\limits_i^n \mathcal{A}_i^2 \right\Vert $.
Then,  given a positive real vector $\mathbf{b} \in \mathbb{R}^p$ with $\tilde{j} \define \arg\min\limits_j \{ b_j \}$ and $t \left( F(x_1,\cdots,x_n) - \mathbb{E}F(x_1,\cdots,x_n )\right)$ satisfing Eq.~\eqref{eq1:lma: Laplace Transform Method Eigentuple Version} for any $t >0$, we have following inequality:
\begin{eqnarray}\label{eq2:thm:Tensor McDiarmid eigentuple}
\mathrm{Pr} \left( \mathbf{d}_{\max}\left( F(x_1,\cdots,x_n) - \mathbb{E}F(x_1,\cdots,x_n )\right)\geq \mathbf{b} \right) \leq mp  e^{-\frac{b_{\tilde{j}}^2}{8 \sigma^2}}.
\end{eqnarray}
\end{restatable}

\subsection{Paper Organization}\label{subsec:Paper Organization}

The rest of this paper is organized as follows. In Section~\ref{sec:Key Results From Part I Paper} , we briefly present those important results from Part I which will be used in later sections. 
Section~\ref{sec:Hermitian T-product Tensors With Random Sequences} utilizes Gaussian and Rademacher series as case studies to explore T-product tensor inequalities. T-product tensor Chernoff bound and its applications are discussed in Section~\ref{sec:Chernoff Bounds for T-Symmetric Tensors}. In Section~\ref{sec:Bernstein Bounds for T-Symmetric Tensors}, T-product tensor Bernstein bound and its applications are provided. Several martingale results based on random T-product  tensors are discussed in Section~\ref{sec:T-Symmetric Tensor Martingales Inequalities}. Concluding remarks are given by Section~\ref{sec:Conclusion}.

\section{Key Results From Part I Paper}\label{sec:Key Results From Part I}\label{sec:Key Results From Part I Paper}

This section will review those important results obtained from Part I paper which will be used at later proofs for references conveneince. All proofs for facts listed in this section can be found at our Part I paper. 

For any tensor $\mathcal{C} \in \mathbb{C}^{m \times n \times p}$, a dilation for the tensor  $\mathcal{C}$, denoted as $\mathfrak{D}(\mathcal{C})$, will be
\begin{eqnarray}\label{eq:dilation T tensor}
\mathfrak{D}(\mathcal{C} )\define  \left[
    \begin{array}{cc}
     \mathcal{O} & \mathcal{C} \\
     \mathcal{C}^{\mathrm{H}} & \mathcal{O} \\
    \end{array}
\right],
\end{eqnarray}
where $\mathfrak{D}(\mathcal{C} ) \in  \mathbb{C}^{(m+n) \times (m+n) \times p}$ and we have $ \left( \mathfrak{D}(\mathcal{C} ) \right)^{\mathrm{H}} =   \mathfrak{D}(\mathcal{C} ) $ (Hermitian T-product tensor after dilation). 

From T-SVD, we have following relation for Hermitian T-product tensor:
\begin{eqnarray}\label{eq:transfer rule}
f(s) \leq g(s)~~ \mbox{for $s \in [a, b]$} \Longrightarrow f(\mathcal{C}) \preceq  g(\mathcal{C}) \mbox{~~when the eigenvalues of $\mathcal{C}$ lie in $[a, b]$.}
\end{eqnarray}
Above Eq.~\eqref{eq:transfer rule} is named as transfer rule.

\begin{corollary}\label{cor:3.3}
Let $\mathcal{A}$ be a fixed Hermitian T-product tensor, and let $\mathcal{X}$ be a random Hermitian T-product tensor, then we have
\begin{eqnarray}
\mathbb{E} \mathrm{Tr} e^{\mathcal{A} + \mathcal{X}} \leq \mathrm{Tr} e^{\mathcal{A} + \log \left( \mathbb{E} e^{\mathcal{X}} \right) }.
\end{eqnarray}
\end{corollary}

\begin{corollary}\label{cor:3_7}
Given a finite sequence of independent Hermitian random tensors $\{ \mathcal{X}_i \} \in \mathbb{C}^{m \times m \times p}$. If there is a function $f: (0, \infty) \rightarrow [0, \infty]$ and a sequence of non-random Hermitian T-product tensors $\{ \mathcal{A}_i \}$ with following condition:
\begin{eqnarray}\label{eq:cond in cor 3_7}
f(t) \mathcal{A}_i \succeq  \log \mathbb{E} e^{t \mathcal{X}_i},~~\mbox{for $t > 0$.}
\end{eqnarray}
Then, for all $\theta \in \mathbb{R}$, we have
\begin{eqnarray}
\mathrm{Pr} \left( \lambda_{\max}\left(\sum\limits_{i=1}^n \mathcal{X}_i \right) \geq \theta \right)
& \leq & mp \inf\limits_{t > 0}\Big\{\exp\left[ - t \theta + f( t ) \lambda_{\max}\left(\sum\limits_{i=1}^n \mathcal{A}_i \right) \right] \Big\}
\end{eqnarray}
\end{corollary}

\begin{corollary}\label{cor:3_7 eigentuple}
Given a finite sequence of independent random Hermitian T-product tensors $\{ \mathcal{X}_i \}$ with dimensions in $\mathbb{C}^{m \times m \times p}$. If there is a function $f: (0, \infty) \rightarrow [0, \infty]$ and a sequence of non-random Hermitian T-product tensors $\{ \mathcal{A}_i \}$ with following condition:
\begin{eqnarray}\label{eq:cond in cor 3_7 eigentuple}
f(t) \mathcal{A}_i \succeq  \log \mathbb{E} e^{t \mathcal{X}_i},~~\mbox{for $t > 0$.}
\end{eqnarray}
Then, for all $\mathbf{b} \in \mathbb{R}^p$ and $\sum\limits_{i=1}^n t \mathcal{X}_i $ satisfing Eq.~\eqref{eq1:lma: Laplace Transform Method Eigentuple Version}, we have
\begin{eqnarray}
\mathrm{Pr} \left( \mathbf{d}_{\max}\left(\sum\limits_{i=1}^n \mathcal{X}_i \right) \geq \mathbf{b} \right)
& \leq & mp \inf\limits_{t > 0} \min\limits_{1 \leq j \leq p} \left\{ \frac{ \exp  \left(  f(t) \lambda_{\max} \left( \sum\limits_{i=1}^n \mathcal{A}_i \right)    \right)    }{ \left( e_{\bigodot}^{t \mathbf{b} }\right)_j } \right\}. 
\end{eqnarray}
\end{corollary}

\begin{corollary}\label{cor:3_9}
Given a finite sequence of independent Hermitian random tensors $\{ \mathcal{X}_i \} \in \mathbb{C}^{m \times m \times p}$. For all $\theta \in \mathbb{R}$, we have
\begin{eqnarray}
\mathrm{Pr} \left( \lambda_{\max}\left(\sum\limits_{i=1}^n \mathcal{X}_i \right) \geq \theta \right)
& \leq &m p  \inf\limits_{t > 0}  \Big\{  \exp\left[ - t \theta +n \log \lambda_{\max} \left( \frac{1}{n} \sum\limits_{i=1}^{n} \mathbb{E} e^{t \mathcal{X}_i} \right)  \right]  \Big\} \nonumber \\
\end{eqnarray}
\end{corollary}

\begin{corollary}\label{cor:3_9 eigentuple}
Given a finite sequence of independent random Hermitian T-product tensors $\{ \mathcal{X}_i \}$ with dimensions in $\mathbb{C}^{m \times m \times p}$, a real vector $\mathbf{b} \in \mathbb{R}^p$ and $\sum\limits_{i=1}^n t \mathcal{X}_i $ satisfing Eq.~\eqref{eq1:lma: Laplace Transform Method Eigentuple Version}, we have
\begin{eqnarray}
\mathrm{Pr} \left( \mathbf{d}_{\max}\left(\sum\limits_{i=1}^n \mathcal{X}_i \right) \geq \mathbf{b} \right)
& \leq & m p \inf\limits_{t > 0} \min\limits_{1 \leq j \leq p} \left\{ \frac{\exp \left( n \log \lambda_{\max}\left( \frac{1}{n} \sum\limits_{i=1}^{n}\mathbb{E}e^{t \mathcal{X}_i } \right) \right)     }{ \left( e_{\bigodot}^{t \mathbf{b} }\right)_j } \right\} 
\end{eqnarray}
\end{corollary}

\begin{lemma}[Laplace Transform Method for T-product Tensors: Eigenvalue Version]\label{lma: Laplace Transform Method Eigenvalue Version}
Let $\mathcal{X}$ be a random Hermitian T-product tensor. For $\theta \in \mathbb{R}$, we have
\begin{eqnarray}
\mathbb{P} (\lambda_{\max}(\mathcal{X}) \geq \theta) \leq \inf_{t > 0} \Big\{ e^{-\theta t} \mathbb{E}\mathrm{Tr} e^{t \mathcal{X}} \Big\}
\end{eqnarray}
\end{lemma}

\begin{lemma}[Laplace Transform Method for T-product Tensors: Eigentuple Version]\label{lma: Laplace Transform Method Eigentuple Version}
Let $\mathcal{X} \in \mathbb{C}^{m \times m \times p}$ be a random T-positive definite (TPD) tensor and an all one vector $\mathbf{1}_p = [1, 1, \cdots, 1]^\mathrm{T} \in \mathbb{C}^{p}$. If $t \mathcal{X}$ satisfies Eq.~\eqref{eq1:lma: Laplace Transform Method Eigentuple Version}, then, for $\mathbf{b} \in \mathbb{R}^p$, we have
\begin{eqnarray}
\mathbb{P} (\mathbf{d}_{\max}(\mathcal{X}) \geq \mathbf{b}) \leq \inf_{t > 0}
 \min\limits_{i} \left\{ \frac{  \mathbb{E} \left( \mathrm{Tr}\left(e^{t \mathcal{X}} \right)\right)   }{ \left(e_{\bigodot}^{ t \mathbf{b}}  \right)_i  }\right\},
\end{eqnarray}
where $\mathbf{d}_{\max}$ is the maximum eigentuple of the TPD tensor $\mathcal{X} $.
\end{lemma}

\begin{theorem}[Golden-Thompson inequality for T-product Tensors]\label{thm:GT Inequality for T-product Tensors}
Given two Hermitian T-product tensors $\mathcal{C}, \mathcal{D} \in \mathbb{C}^{m \ times m \times p}$, we have 
\begin{eqnarray}
\mathrm{Tr} \left( \exp(\mathcal{C} + \mathcal{D}) \right) \leq \mathrm{Tr} \left( \exp \left( \mathcal{C} \right)  \star \exp \left( \mathcal{D} \right) \right)
\end{eqnarray}
\end{theorem}

\section{Hermitian T-product Tensors With Random Sequences}\label{sec:Hermitian T-product Tensors With Random Sequences}

A Hermitian T-product tensor Gaussian series is one of the simplest cases of a sum of independent random Hermitian T-product tensors. For scalers, a Gaussian series with real coefficients satisfies a normal-type tail bound where the variance is controlled by the sum of squares coefficients. The first Section~\ref{sec:Tensor with Gaussian and Rademacher Random Series} is to extend this context to Hermitian T-product tensors. In Section~\ref{sec:A Gaussian Tensor with Nonuniform Variances}, we will apply results from Section~\ref{sec:Tensor with Gaussian and Rademacher Random Series} to consider Gaussian Hermitian T-product tensor with nonuniform variances. Finally, we will provide the lower and upper bounds of random Hermitian T-product tensor expectation in Section~\ref{sec:Lower and Upper Bounds of Tensor Expectation}.

\subsection{Hermitian T-product Tensors with Gaussian and Rademacher Random Series}\label{sec:Tensor with Gaussian and Rademacher Random Series}

We begin with a lemma about moment-generating functions of Rademacher and Gaussian normal random variables. 

\begin{lemma}\label{lma:mgf of Rademacher and normal rvs.}
Suppose that the tensor $\mathcal{A} \in \mathbb{C}^{m \times m \times p}$ is Hermitian T-product tensor.  Given a Gaussian normal random variable $\alpha$ and a Rademacher random variable $\beta$, then, we have
\begin{eqnarray}
\mathbb{E}e^{\alpha t \mathcal{A}} =  e^{t^2 \mathcal{A}^2 / 2}  \mbox{~~and~~} e^{t^2 \mathcal{A}^2 / 2} \succeq \mathbb{E}e^{\beta t \mathcal{A}},
\end{eqnarray}   
where $t \in \mathbb{R}$.
\end{lemma}
\textbf{Proof:}
For the standard normal random variable, because we have
\begin{eqnarray}
\mathbb{E} (\alpha^{2i}) = \frac{ (2i)!}{i ! 2^i}  \mbox{~~and~~} \mathbb{E} (\alpha^{2i + 1}) = 0, 
\end{eqnarray}
where $i = 0,1,2,\cdots$; then
\begin{eqnarray}
\mathbb{E}e^{\alpha t \mathcal{A}} &=& \mathcal{I}_{mmp} + \sum\limits_{i=1}^{\infty} \frac{  \mathbb{E}(\alpha^{2i}) (t \mathcal{A})^{2i}}{ (2i)! } \nonumber \\
&=& \mathcal{I}_{mmp} + \sum\limits_{i=1}^{\infty} \frac{ (t^2 \mathcal{A}^2/2)^i}{ i! } = e^{t^2 \mathcal{A}^2 /2}. 
\end{eqnarray}

For the Rademacher random variable, we have 
\begin{eqnarray}
\mathbb{E} e^{\beta  t \mathcal{A}} = \cosh (t \mathcal{A}) \preceq e^{t^2 \mathcal{A}^2/2}.
\end{eqnarray}
Therefore, this Lemma is proved.
$\hfill \Box$

We are ready to present the main theorem of this section about Hermitian T-product tensors with Gaussian and Rademacher series. The eigenvalue version is provided first by Theorem~\ref{thm:TensorGaussianNormalSeries eigenvalue}.

\TensorGaussianNormalSeriesEigenvalue*

%
\textbf{Proof:}
Given a finite sequence of independent Gaussian or Rademacher random variables $\{ \alpha_i \}$, from Lemma~\ref{lma:mgf of Rademacher and normal rvs.}, we have 
\begin{eqnarray}
e^{\frac{t^2 \mathcal{A}_i^2}{2}} \succeq \mathbb{E}e^{\alpha_i t \mathcal{A}_i}.
\end{eqnarray}
From the definition in Eq.~\eqref{eq:4_5} and Corollary~\ref{cor:3_7}, we have 
\begin{eqnarray}\label{eq:4_9}
\mathrm{Pr}\left( \lambda_{\max} \left( \sum\limits_{i=1}^{n}\alpha_i \mathcal{A}_i \right) \geq \theta \right) \leq  mp \inf\limits_{t > 0}\Big\{ e^{- t \theta + \frac{t^2 \sigma^2}{2} } \Big\} = mp e^{- \frac{\theta^2}{2 \sigma^2}}.
\end{eqnarray}
This establishes  Eq.~\eqref{eq:4_3}. For  Eq.~\eqref{eq:4_4}, we have to apply following facts about the symmetric distribution of Gaussian and Rademacher random variables to obtain  
\begin{eqnarray}
\mathrm{Pr}\left( \lambda_{\max} \left( \sum\limits_{i=1}^{n} (-\alpha_i) \mathcal{A}_i \right) \geq \theta \right) =
\mathrm{Pr}\left(  - \lambda_{\min} \left( \sum\limits_{i=1}^{n}\alpha_i \mathcal{A}_i \right) \geq \theta \right) \leq mp e^{- \frac{\theta^2}{2 \sigma^2}}.
\end{eqnarray}
Then, we obtain  Eq.~\eqref{eq:4_4} as follows: 
\begin{eqnarray}
\mathrm{Pr}\left( \left\Vert \sum\limits_{i=1}^{n}\alpha_i \mathcal{A}_i \right\Vert \geq \theta \right)
&=& 2 \mathrm{Pr}\left( \lambda_{\max} \left( \sum\limits_{i=1}^{n}\alpha_i \mathcal{A}_i \right) \geq \theta \right)
\nonumber \\  
& \leq &  2 mp e^{-\frac{\theta^2 }{2 \sigma^2}}.
\end{eqnarray}
$\hfill \Box$

From the Hermitian dilation definition provided by Eq.~\eqref{eq:dilation T tensor}, we can extend Theorem~\ref{thm:TensorGaussianNormalSeries eigenvalue} from square Hermitian tensor to rectangular tensor by the following corollary.
\begin{corollary}[Rectangular Tensor with Gaussian and Rademacher Series Eigenvalue Version]\label{cor:Rectangular Tensor with Gaussian and Rademacher Series eigenvalue}
Given a finite sequence $\mathcal{A}_i \in \mathbb{C}^{m \times n \times p}$ be a finite sequence of indepedent standard normal random variables. We define 
\begin{eqnarray}\
\sigma^2 &\define& \max \Bigg\{ \left\Vert \sum\limits_{i=1}^n \mathcal{A}_i \star \mathcal{A}^{\mathrm{H}}_i \right\Vert,  \left\Vert \sum\limits_{i=1}^n \mathcal{A}^{\mathrm{H}}_i \star \mathcal{A}_i \right\Vert\Bigg\}.
\end{eqnarray}
then, for all $\theta \geq 0$, we have 
\begin{eqnarray}
\mathrm{Pr}\left( \left\Vert \sum\limits_{i=1}^n \alpha_i \mathcal{A}_i \right\Vert \geq \theta \right) \leq  (m+n) p  e^{-\frac{\theta^2 }{2 \sigma^2}}.
\end{eqnarray}
This corollary is also valid for a finite sequence of independent Rademacher random variables $\{ \alpha_i \}$.
\end{corollary}
\textbf{Proof:}
Let $\{ \alpha_i \}$ be a finite sequence of independent Gaussian or Rademacher random variables. Consider a finite sequence of random Hermitian T-product tensors $\{ \alpha_i \mathfrak{D}(\mathcal{A}_i) \}$ with dimensions $\mathbb{C}^{ (m+n) \times (m+n) \times p}$ and the fact that the largest eigenvalue of $\mathfrak{D}(\mathcal{A}_i )$ will be the same with the largest singular of $\mathcal{A}_i$,  we have 
\begin{eqnarray}\label{eq1:cor:Rectangular Tensor Gaussian with Normal Series}
\left\Vert \sum\limits_i^n \alpha_i \mathcal{A}_i \right\Vert &=& \lambda_{\max} \left( \mathbb{D} \left(  \sum\limits_{i=1}^n\alpha_i \mathcal{A}_i \right)\right) = \lambda_{\max} \left(\sum\limits_{i=1}^n \alpha_i  \mathbb{D} \left(\mathcal{A}_i \right)\right).
\end{eqnarray}
Due to the following singular value relation
\begin{eqnarray}\label{eq2:cor:Rectangular Tensor Gaussian with Normal Series}
\sigma^2 &=& \left \Vert  \sum\limits_i^n \mathbb{D}( \mathcal{A}_i)^2\right\Vert
=
\left\Vert
\begin{bmatrix}
\sum\limits_{i=1}^n \mathcal{A}_i \star \mathcal{A}_i^{\mathrm{H}} & \mathcal{O}   \\
\mathcal{O} &\sum\limits_{i=1}^n \mathcal{A}^{\mathrm{H}}_i \star \mathcal{A}_i \\
\end{bmatrix}\right\Vert \nonumber \\
&=& \max \Bigg\{ \left\Vert\sum\limits_{i=1}^n  \mathcal{A}_i \star \mathcal{A}^{\mathrm{H}}_i \right\Vert,  \left\Vert \sum\limits_i^{n} \mathcal{A}^{\mathrm{H}}_i \star \mathcal{A}_i \right\Vert\Bigg\}.
\end{eqnarray}
From Eqs.~\eqref{eq1:cor:Rectangular Tensor Gaussian with Normal Series}, and Theorem~\ref{thm:TensorGaussianNormalSeries eigenvalue}, this corollary is proved.
$\hfill \Box$

The eigentuple version for Theorem~\ref{thm:TensorGaussianNormalSeries eigenvalue} is provided by the folloiwing Theorem~\ref{thm:TensorGaussianNormalSeries eigentuple}.

\TensorGaussianNormalSeriesEigentuple*

%
\textbf{Proof:}
Given a finite sequence of independent Gaussian or Rademacher random variables $\{ \alpha_i \}$, from Lemma~\ref{lma:mgf of Rademacher and normal rvs.}, we have 
\begin{eqnarray}
e^{\frac{t^2 \mathcal{A}_i^2}{2}} \succeq \mathbb{E}e^{\alpha_i t \mathcal{A}_i}.
\end{eqnarray}

If $\tilde{j}$ is determined as:
\begin{eqnarray}\label{eq:tilde i def}
\tilde{j} &=& \arg \min\limits_j \left\{ b_j \right\},
\end{eqnarray}
where $b_j$ are entries of the vector $\mathbf{b}$. Then, we have  
\begin{eqnarray}
 \min_{1 \leq j \leq p}\left\{\frac{ \exp \left(\frac{t^2}{2} \lambda_{\max} \left( \sum\limits_{i=1}^{n}  \mathcal{A}_i^2   \right) \right) }{   \left( e_{\bigodot}^{t \mathbf{b}} \right)_{j} } \right\}
& \leq & 
\exp \left(  - t b_{\tilde{j}} +  \frac{t^2}{2} \lambda_{\max} \left( \sum\limits_{i=1}^{n}  \mathcal{A}_i^2   \right) \right). 
\end{eqnarray}

From the definition in Eq.~\eqref{eq:4_5 eigentuple} and Corollary~\ref{cor:3_7 eigentuple}, we have 
\begin{eqnarray}\label{eq:4_9}
\mathrm{Pr}\left( \mathbf{d}_{\max} \left( \sum\limits_{i=1}^{n}\alpha_i \mathcal{A}_i \right) \geq \mathbf{b} \right) \leq  mp \inf\limits_{t > 0}\Big\{ e^{- t b_{\tilde{j}} + \frac{t^2 \sigma^2}{2} } \Big\} = mp e^{- \frac{b_{\tilde{j}}^2}{2 \sigma^2}}. 
\end{eqnarray}

For  Eq.~\eqref{eq:4_4 eigentuple}, because Gaussian and Rademacher random variables are symmetric, we have 
\begin{eqnarray}
\mathrm{Pr}\left( \mathbf{d}_{\max} \left( \sum\limits_{i=1}^{n} (-\alpha_i) \mathcal{A}_i \right) \geq \mathbf{b} \right) =
\mathrm{Pr}\left(  - \mathbf{d}_{\min} \left( \sum\limits_{i=1}^{n}\alpha_i \mathcal{A}_i \right) \geq \mathbf{b} \right) \leq mp e^{- \frac{b_{\tilde{j}}^2}{2 \sigma^2}}.
\end{eqnarray}
Then, we obtain Eq.~\eqref{eq:4_4 eigentuple} as follows: 
\begin{eqnarray}
\mathrm{Pr}\left( \left\Vert \sum\limits_{i=1}^{n}\alpha_i \mathcal{A}_i \right\Vert_{\mbox{\tiny{vec}}} \geq \mathbf{b} \right)
&=& 2 \mathrm{Pr}\left( \mathbf{d}_{\max} \left( \sum\limits_{i=1}^{n}\alpha_i \mathcal{A}_i \right) \geq \mathbf{b} \right)
\nonumber \\  
& \leq &  2 mp e^{- \frac{b_{\tilde{j}}^2}{2 \sigma^2}}.
\end{eqnarray}
$\hfill \Box$

From the Hermitian dilation definition provided by Eq.~\eqref{eq:dilation T tensor}, we can extend Theorem~\ref{thm:TensorGaussianNormalSeries eigentuple} from square Hermitian tensor to rectangular tensor by the following corollary.
\begin{corollary}[Rectangular Tensor with Gaussian and Rademacher Series Eigentuple Version]\label{cor:Rectangular Tensor with Gaussian and Rademacher Series eigentuple}
Given a finite sequence $\mathcal{A}_i \in \mathbb{C}^{m \times n \times p}$ be a finite sequence of indepedent standard normal random variables. We define 
\begin{eqnarray}\label{eq:sigma 2 def cor:Rectangular Tensor with Gaussian and Rademacher Series eigentuple}
\sigma^2 &\define& \max \Bigg\{ \left\Vert \sum\limits_{i=1}^n \mathcal{A}_i \star \mathcal{A}^{\mathrm{H}}_i \right\Vert,  \left\Vert \sum\limits_{i=1}^n \mathcal{A}^{\mathrm{H}}_i \star \mathcal{A}_i \right\Vert\Bigg\}.
\end{eqnarray}
then, for all $\mathbf{b} \geq \mathbf{0}$ and $\sum\limits_{i=1}^{n}t  \alpha_i \mathcal{A}_i$ satisfying Eq.~\eqref{eq1:lma: Laplace Transform Method Eigentuple Version} for $t >0$, we have 
\begin{eqnarray}
\mathrm{Pr}\left(  \left\Vert \sum\limits_{i=1}^{n} \alpha_i \mathcal{A}_i \right\Vert_{\mbox{\tiny{vec}}} \geq \mathbf{b} \right) \leq  (m+n) p  e^{-\frac{b_{\tilde{j}}^2 }{2 \sigma^2}},
\end{eqnarray}
where $\tilde{j}$ is defined by Eq.~\eqref{eq:tilde i def}.

This corollary is also valid for a finite sequence of independent Rademacher random variables $\{ \alpha_i \}$.
\end{corollary}
\textbf{Proof:}
Let $\{ \alpha_i \}$ be a finite sequence of independent Gaussian or Rademacher random variables. Consider a finite sequence of random Hermitian T-product tensors $\{ \alpha_i \mathfrak{D}(\mathcal{A}_i) \}$ with dimensions $\mathbb{C}^{ (m+n) \times (m+n) \times p}$ and the fact that the largest eigentuple of $\mathfrak{D}(\mathcal{A}_i )$ will be the same with the largest eigentuple of $\mathcal{A}_i$,  we have 
\begin{eqnarray}\label{eq1:cor:Rectangular Tensor with Gaussian and Rademacher Series  eigentuple}
\left\Vert  \sum\limits_{i=1}^{n} \alpha_i \mathcal{A}_i \right\Vert_{\mbox{\tiny{vec}}} &=& \mathbf{d}_{\max} \left( \mathbb{D} \left(  \sum\limits_{i=1}^n\alpha_i \mathcal{A}_i \right)\right)=\mathbf{d}_{\max} \left(\sum\limits_{i=1}^n \alpha_i  \mathbb{D} \left(\mathcal{A}_i \right)\right).
\end{eqnarray}
Due to the following singular value relation
\begin{eqnarray}\label{eq2:cor:Rectangular Tensor Gaussian with Normal Series}
\sigma^2 &=& \left \Vert  \sum\limits_i^n \mathbb{D}( \mathcal{A}_i)^2\right\Vert
=
\left\Vert
\begin{bmatrix}
\sum\limits_{i=1}^n \mathcal{A}_i \star \mathcal{A}_i^{\mathrm{H}} & \mathcal{O}   \\
\mathcal{O} &\sum\limits_{i=1}^n \mathcal{A}^{\mathrm{H}}_i \star \mathcal{A}_i \\
\end{bmatrix}\right\Vert \nonumber \\
&=& \max \Bigg\{ \left\Vert\sum\limits_{i=1}^n  \mathcal{A}_i \star \mathcal{A}^{\mathrm{H}}_i \right\Vert,  \left\Vert \sum\limits_i^{n} \mathcal{A}^{\mathrm{H}}_i \star \mathcal{A}_i \right\Vert\Bigg\}.
\end{eqnarray}
From Eq.~\eqref{eq1:cor:Rectangular Tensor with Gaussian and Rademacher Series  eigentuple} and Theorem~\ref{thm:TensorGaussianNormalSeries eigentuple}, this corollary is proved.
$\hfill \Box$

\subsection{A Gaussian Tensor with Nonuniform Variances}\label{sec:A Gaussian Tensor with Nonuniform Variances}

In this section, we will apply results obtained from the previous section to consider Gaussian tensor with nonuniform variances among random entries. 
\begin{corollary}\label{cor:Gaussian Tensor with Nonuniform Var eigenvalue}
Given a tensor $\mathcal{A} \in \mathbb{C}^{m \times n \times p}$ and a random tensor $\mathcal{X} \in \mathbb{C}^{m \times n \times p}$ whose entries are independent standard Gaussian normal random variables. Let $\circ$ represent the Hadamard product (entrywise) between two T-product tensors with the same dimensions. Then, we have
\begin{eqnarray}
\mathrm{Pr} \left( \left\Vert \mathcal{X} \circ \mathcal{A} \right\Vert \geq \theta \right) 
\leq (m+n)p e^{-\frac{\theta^2 }{2 \sigma^2}}, 
\end{eqnarray}
where 
\begin{eqnarray}
\sigma^2 &=& \max \left\{ \sum\limits_{j=1}^{n} \left\vert a_{1, j, 1}\right\vert^2, \sum\limits_{j=1}^{n}\left\vert a_{2, j, 1}\right\vert^2, \cdots, \sum\limits_{j=1}^{n} \left\vert a_{m, j, 1}\right\vert^2, \right. \nonumber \\
&  & 
\left.
~~~~~~~~~~~ \sum\limits_{i=1}^{m} \left\vert a_{i, 1, 1}\right\vert^2, \sum\limits_{i=1}^{m}\left\vert a_{i, 2, 1}\right\vert^2, \cdots, \sum\limits_{i=1}^{m} \left\vert a_{i, n, 1}\right\vert^2  \right\}
\end{eqnarray}
where $a_{i, j, k}$ are entries of the tensor $\mathcal{A}$.
\end{corollary}
\textbf{Proof:}
Since we can decompose the tensor $\mathcal{X} \circ \mathcal{A}$ as:
\begin{eqnarray}
\mathcal{X} \circ \mathcal{A}  &=& \sum\limits_{i=j=k=1}^{m, n, p} x_{i, j, k} a_{i, j, k} \mathcal{E}_{i, j, k},  
\end{eqnarray}
where $\mathcal{E}_{i, j, k} \in \mathbb{C}^{m \times n \times p}$ is the tensor with all zero entries except unity at the position $i, j, k$; then, we have
\begin{eqnarray}
\sum\limits_{i=j=k=1}^{m, n, p} \left(a_{i, j, k} \mathcal{E}_{i, j, k} \right) \star \left(a_{i, j, k} \mathcal{E}_{i, j, k} \right)^{\mathrm{H}} 
&=& \sum\limits_{i=k=1}^{m, p}\left( \sum\limits_{j=1}^n \left\vert a_{i, j, k} \right\vert^2 \right)     \mathcal{E}_{i, i, 1}  \nonumber \\
&=& \mbox{fdiag} \left( \sum\limits_{j=1}^{n} \left\vert a_{1, j, 1}\right\vert^2, \sum\limits_{j=1}^{n}\left\vert a_{2, j, 1}\right\vert^2, \cdots, \sum\limits_{j=1}^{n} \left\vert a_{m, j, 1}\right\vert^2\right),
\end{eqnarray}
where $\mbox{fdiag}$ is the tensor with dimensions in $\mathbb{C}^{m \times m \times p}$ such that the frontal diagonal matrix is a diagonal matrix and zero matrices at the other matrices parallel to the frontal matrix. Similarly, we also have 
\begin{eqnarray}
\sum\limits_{i=j=k=1}^{m, n, p} \left(a_{i, j, k} \mathcal{E}_{i, j, k} \right)^{\mathrm{H}}  \star \left(a_{i, j, k} \mathcal{E}_{i, j, k} \right)
&=& \sum\limits_{j=k=1}^{n, p}\left( \sum\limits_{i=1}^m \left\vert a_{i, j, k} \right\vert^2 \right)     \mathcal{E}_{j, j, 1}  \nonumber \\
&=& \mbox{fdiag} \left( \sum\limits_{i=1}^{m} \left\vert a_{i, 1, 1}\right\vert^2, \sum\limits_{i=1}^{m}\left\vert a_{i, 2, 1}\right\vert^2, \cdots, \sum\limits_{i=1}^m \left\vert a_{i, n, 1} \right\vert^2 \right),
\end{eqnarray}

Therefore, we have 
\begin{eqnarray}
\sigma^2 &=& \max\left\{ \lambda_{\max} \left( \mbox{fdiag} \left( \sum\limits_{j=1}^{n} \left\vert a_{1, j, 1}\right\vert^2, \sum\limits_{j=1}^{n}\left\vert a_{2, j, 1}\right\vert^2, \cdots, \sum\limits_{j=1}^{n} \left\vert a_{m, j, 1}\right\vert^2\right) \right), \right. \nonumber \\
& & \left.
~~~~~~~~~~~ \lambda_{\max} \left(\mbox{fdiag} \left( \sum\limits_{i=1}^{m} \left\vert a_{i, 1, 1}\right\vert^2, \sum\limits_{i=1}^{m}\left\vert a_{i, 2, 1}\right\vert^2, \cdots, \sum\limits_{i=1}^m \left\vert a_{i, n, 1} \right\vert^2 \right) \right) \right\}\nonumber \\
&=&  \max \left\{ \sum\limits_{j=1}^{n} \left\vert a_{1, j, 1}\right\vert^2, \sum\limits_{j=1}^{n}\left\vert a_{2, j, 1}\right\vert^2, \cdots, \sum\limits_{j=1}^{n} \left\vert a_{m, j, 1}\right\vert^2, \right. \nonumber \\
&  & 
\left.
~~~~~~~~~~~ \sum\limits_{i=1}^{m} \left\vert a_{i, 1, 1}\right\vert^2, \sum\limits_{i=1}^{m}\left\vert a_{i, 2, 1}\right\vert^2, \cdots, \sum\limits_{i=1}^{m} \left\vert a_{i, n, 1}\right\vert^2  \right\}
\end{eqnarray}
Finally, from Corollary~\ref{cor:Rectangular Tensor with Gaussian and Rademacher Series eigenvalue}, this Corollary is proved.
$\hfill \Box$

Following corollary is the eigentuple version for Corollary~\ref{cor:Gaussian Tensor with Nonuniform Var eigenvalue}

\begin{corollary}\label{cor:Gaussian Tensor with Nonuniform Var eigentuple}
Given a tensor $\mathcal{A} \in \mathbb{C}^{m \times n \times p}$ and a random tensor $\mathcal{X} \in \mathbb{C}^{m \times n \times p}$ whose entries are independent standard Gaussian normal random variables. Let $\circ$ be used to represent a Hadamard product (entrywise) between two T-product tensors with the same dimensions. Then, for all $\mathbf{b} \geq \mathbf{0}$ with $\tilde{j}$ defined by Eq.~\eqref{eq:tilde i def}, and $ t \mathcal{X} \circ \mathcal{A} $ satisfying Eq.~\eqref{eq1:lma: Laplace Transform Method Eigentuple Version} for $t >0$, we have
\begin{eqnarray}
\mathrm{Pr} \left( \left\Vert \mathcal{X} \circ \mathcal{A} \right\Vert_{\mbox{\tiny{vec}}} \geq \mathbf{b} \right) 
\leq (m+n)p e^{-\frac{b_{\tilde{j}}^2 }{2 \sigma^2}}, 
\end{eqnarray}
where 
\begin{eqnarray}
\sigma^2 &=& \max \left\{ \sum\limits_{j=1}^{n} \left\vert a_{1, j, 1}\right\vert^2, \sum\limits_{j=1}^{n}\left\vert a_{2, j, 1}\right\vert^2, \cdots, \sum\limits_{j=1}^{n} \left\vert a_{m, j, 1}\right\vert^2, \right. \nonumber \\
&  & 
\left.
~~~~~~~~~~~ \sum\limits_{i=1}^{m} \left\vert a_{i, 1, 1}\right\vert^2, \sum\limits_{i=1}^{m}\left\vert a_{i, 2, 1}\right\vert^2, \cdots, \sum\limits_{i=1}^{m} \left\vert a_{i, n, 1}\right\vert^2  \right\}.
\end{eqnarray}
The terms $a_{i, j, k}$ are entries of the tensor $\mathcal{A}$.

\end{corollary}
\textbf{Proof:}
From the proof from Corollary~\ref{cor:Gaussian Tensor with Nonuniform Var eigenvalue}, we also have 
\begin{eqnarray}
\sigma^2 &=& \max\left\{ \lambda_{\max} \left( \mbox{fdiag} \left( \sum\limits_{j=1}^{n} \left\vert a_{1, j, 1}\right\vert^2, \sum\limits_{j=1}^{n}\left\vert a_{2, j, 1}\right\vert^2, \cdots, \sum\limits_{j=1}^{n} \left\vert a_{m, j, 1}\right\vert^2\right) \right), \right. \nonumber \\
& & \left.
~~~~~~~~~~~ \lambda_{\max} \left(\mbox{fdiag} \left( \sum\limits_{i=1}^{m} \left\vert a_{i, 1, 1}\right\vert^2, \sum\limits_{i=1}^{m}\left\vert a_{i, 2, 1}\right\vert^2, \cdots, \sum\limits_{i=1}^m \left\vert a_{i, n, 1} \right\vert^2 \right) \right) \right\}\nonumber \\
&=&  \max \left\{ \sum\limits_{j=1}^{n} \left\vert a_{1, j, 1}\right\vert^2, \sum\limits_{j=1}^{n}\left\vert a_{2, j, 1}\right\vert^2, \cdots, \sum\limits_{j=1}^{n} \left\vert a_{m, j, 1}\right\vert^2, \right. \nonumber \\
&  & 
\left.
~~~~~~~~~~~ \sum\limits_{i=1}^{m} \left\vert a_{i, 1, 1}\right\vert^2, \sum\limits_{i=1}^{m}\left\vert a_{i, 2, 1}\right\vert^2, \cdots, \sum\limits_{i=1}^{m} \left\vert a_{i, n, 1}\right\vert^2  \right\}
\end{eqnarray}
Finally, from Corollary~\ref{cor:Rectangular Tensor with Gaussian and Rademacher Series eigentuple}, this Corollary is proved.
$\hfill \Box$

\subsection{Lower and Upper Bounds of Spectral Norm Expectation}\label{sec:Lower and Upper Bounds of Tensor Expectation}
Given a finite sequence $\mathcal{A}_i \in \mathbb{C}^{m \times m \times p}$, and let $\{ \alpha_i \}$ be a finite sequence of indepedent standard normal variables. We define following random tensor 
\begin{eqnarray}
\mathcal{X} = \sum\limits_{i=1}^n \alpha_i \mathcal{A}_i.
\end{eqnarray}
From Theorem~\ref{thm:TensorGaussianNormalSeries eigenvalue}, we have 
\begin{eqnarray}\label{eq:upper bound for exp tensor}
\mathbb{E}\left( \left\Vert \mathcal{X} \right\Vert^2 \right)& = & \int_0^{\infty} \mathrm{Pr}\left(
\left\Vert \mathcal{X} \right\Vert > \sqrt{t} \right) dt \leq  \int_0^{\infty} 2 mp e^{-\frac{t}{2 \sigma^2}} dt = 4mp\sigma^2
\end{eqnarray}
where $\sigma^2 = \left\Vert \sum\limits_{i=1}^n \mathcal{A}_i^2 \right\Vert$. On the other hand, from Jensen's inequality, we have 
\begin{eqnarray}\label{eq:lower bound for exp tensor}
\mathbb{E}\left( \left\Vert \mathcal{X} \right\Vert^2 \right) = \mathbb{E} \left\Vert \mathcal{X}^2 \right\Vert  \geq \left\Vert \mathbb{E}(\mathcal{X}^2) \right\Vert =  \left\Vert \sum\limits_{i=1}^n \mathcal{A}^2_i\right\Vert = \sigma^2.
\end{eqnarray}
From both Eqs.~\eqref{eq:upper bound for exp tensor} and~\eqref{eq:lower bound for exp tensor},
we have following relation:
\begin{eqnarray}
c \sigma \leq \mathbb{E} \left\Vert \mathcal{X} \right\Vert \leq 2\sigma \sqrt{mp}
\end{eqnarray}
This shows that the tensor variance parameter $\sigma^2$ controls the expected norm $ \mathbb{E} \left\Vert \mathcal{X} \right\Vert$ with square root of logarithmic function for the tensor dimensions. 

\section{Chernoff Bounds for T-product Tensors}\label{sec:Chernoff Bounds for T-Symmetric Tensors}

In this section, we will extend Chernoff bounds of random variables to random T-product tensors. 

\subsection{T-product Tensor Chernoff Bounds Derivations}\label{sec:T-product Tensor Chernoff Bounds Derivations}

We begin to present a lemma about the semidefinite relation for the tensor moment-generating function of a random TPSD T-product tensor. 

\begin{lemma}\label{lma:Chernoff MGF}
Given a random TPSD T-product tensor with $\lambda_{\max}(\mathcal{X})\leq 1$, then, for any $t \in \mathbb{R}$, we have 
\begin{eqnarray}
 \mathcal{I}+ (e^t - 1) \mathbb{E} \mathcal{X} \succeq \mathbb{E} e^{t \mathcal{X}}.
\end{eqnarray}
\end{lemma}
\textbf{Proof:}
Consider a convex function $f(x) = e^{t x}$, we have  
\begin{eqnarray}
1 + (e^t - 1)x \geq f(x),
\end{eqnarray}
where $x \in [0,1]$. Since the eigenvalues of the random tensor $\mathcal{X}$ lie in the interval $[0, 1]$, from Eq.~\eqref{eq:transfer rule}, we obtain
\begin{eqnarray}
\mathcal{I} + (e^t - 1)\mathcal{X} \succeq e^{t \mathcal{X}}.
\end{eqnarray}
Then, this Lemma is proved by taking the expectation with respect to the random T-product tensor $\mathcal{X}$.
$\hfill \Box$

Given two real values $c, d \in [0, 1]$, we define \emph{binary information divergence} of $c$ and $d$, expressed by $\mathfrak{D}(c || d)$, as 
\begin{eqnarray}\label{eq:def of binary info div}
\mathfrak{D}(c || d) \define c \log \frac{c}{d} + (1-c) \frac{1-c}{1-d}.
\end{eqnarray}    

We are ready to present T-product tensor Chernoff inequality by~\cref{thm:TensorChernoffBoundI}. 

\TensorChernoffBoundIEigenvalue*

%
%
\textbf{Proof:}
From Lemma~\ref{lma:Chernoff MGF}, we have 
\begin{eqnarray}
\mathcal{I} + f(t) \mathbb{E} \mathcal{X}_i \succeq \mathbb{E} e^{t \mathcal{X}_i},
\end{eqnarray}
where $f(t) \define e^t -1$ for $t > 0$.
By applying Corollary~\ref{cor:3_9}, we obtain
\begin{eqnarray}\label{eq1:Chernoff I Upper Bound proof}
\mathrm{Pr} \left( \lambda_{\max}\left( \sum\limits_{i=1}^{n} \mathcal{X}_i \right) \geq \alpha \right) &\leq& mp \inf\limits_{t > 0} \exp\left( - t \alpha + n \log\lambda_{\max}\left( \frac{1}{n} \sum\limits_{i=1}^{n} \left( \mathcal{I} + f(t) \mathbb{E} \mathcal{X}_i \right) \right) \right) \nonumber \\
&= &  mp  \inf\limits_{t > 0} \exp\left( - t \alpha+ n \log\lambda_{\max}\left(  \mathcal{I} + f(t)  \frac{1}{n} \sum\limits_{i=1}^{n} \mathbb{E} \mathcal{X}_i \right) \right) \nonumber \\ 
&=& mp  \inf\limits_{t > 0} \exp\left( - t \alpha + n \log \left(1 + f(t) \overline{\mu}_{\max} \right) \right) .
\end{eqnarray}
The last equality follows from the definition of $\overline{\mu}_{\max}$ and the eigenvalue map properties. When the value $t$ at the right-hand side of  Eq.~\eqref{eq1:Chernoff I Upper Bound proof} is 
\begin{eqnarray}\label{eq2:Chernoff I Upper Bound proof}
t = \log \frac{\alpha}{n - \alpha} - \log \frac{\overline{\mu}_{\max}}{1 - \overline{\mu}_{\max}},
\end{eqnarray}
we can achieve the tightest upper bound at Eq.~\eqref{eq1:Chernoff I Upper Bound proof}. By substituting the value $t$ in Eq.~\eqref{eq2:Chernoff I Upper Bound proof} into Eq.~\eqref{eq1:Chernoff I Upper Bound proof} and change the variable $\alpha \rightarrow n \theta$, Eq.~\eqref{eq:Chernoff I Upper Bound} is proved. The next goal is to prove Eq.~\eqref{eq:Chernoff I Lower Bound}.

If we apply Lemma~\ref{lma:Chernoff MGF} to the sequence $\{- \mathcal{X}_i \}$, we have 
\begin{eqnarray}
\mathcal{I} - g(t) \mathbb{E} \mathcal{X}_i \succeq \mathbb{E} e^{t (-\mathcal{X}_i)},
\end{eqnarray}
where $g(t) \define 1 - e^t$ for $t > 0$.
By applying Corollary~\ref{cor:3_9} again, we obtain
\begin{eqnarray}\label{eq1:Chernoff I Lower Bound proof}
\mathrm{Pr} \left( \lambda_{\min}\left( \sum\limits_{i=1}^{n} \mathcal{X}_i \right) \leq \alpha \right) &=& \mathrm{Pr} \left( \lambda_{\max}\left( \sum\limits_{i=1}^{n} \left( - \mathcal{X}_i \right)\right) \geq \alpha \right) 
\nonumber \\
&\leq& mp  \inf\limits_{t > 0} \exp\left( t \alpha + n \log\lambda_{\max}\left( \frac{1}{n} \sum\limits_{i=1}^{n} \left( \mathcal{I} - g(t) \mathbb{E} \mathcal{X}_i \right) \right) \right) \nonumber \\
&=_1 &  mp   \inf\limits_{t > 0} \exp\left( t \alpha + n \log \left( 1 -f(t)  \lambda_{\min}\left(    \frac{1}{n} \sum\limits_{i=1}^{n} \mathbb{E} \mathcal{X}_i \right) \right) \right) \nonumber \\ 
&=& mp  \inf\limits_{t > 0} \exp\left( t \alpha + n \log \left(1 - g(t) \overline{\mu}_{\min} \right) \right),
\end{eqnarray}
where we apply the relation $\lambda_{\min}( -\frac{1}{n} \sum\limits_{i=1}^{n} \mathbb{E} \mathcal{X}_i  ) = -\lambda_{\max} (\frac{1}{n} \sum\limits_{i=1}^{n} \mathbb{E} \mathcal{X}_i )$ at the equality $=_1$. When the value $t$ at the right-hand side of  Eq.~\eqref{eq1:Chernoff I Lower Bound proof} is 
\begin{eqnarray}\label{eq2:Chernoff I Lower Bound proof}
t = \log \frac{\overline{\mu}_{\min}}{1 - \overline{\mu}_{\min}} -  \log \frac{\alpha}{n - \alpha},
\end{eqnarray}
we can achieve the tightest upper bound at Eq.~\eqref{eq1:Chernoff I Lower Bound proof}. By substituting the value $t$ in Eq.~\eqref{eq2:Chernoff I Lower Bound proof} into Eq.~\eqref{eq1:Chernoff I Lower Bound proof} and change the variable $\alpha \rightarrow n \theta$, Eq.~\eqref{eq:Chernoff I Lower Bound} is proved also.
$\hfill \Box$

The tensor Chernoff bounds discussed at Theorem~\ref{thm:TensorChernoffBoundI} is not related to $\mu_{\max}$ and $\mu_{\min}$ directly. Following theorem is another version of tensor Chernoff bounds to associate the probability range in terms of $\mu_{\max}$ and $\mu_{\min}$ directly and this format of tensor Chernoff bounds is easier to be applied. 

\TensorChernoffBoundIIEigenvalue*

%
\textbf{Proof:}
Without loss of generality, we can assume $T=1$ in our proof. From  Eq.~\eqref{eq1:Chernoff I Upper Bound proof} and the inequality $\log(1 + x) \leq x$ for $x > -1$, we have 
\begin{eqnarray}\label{eq1:Chernoff II Upper Bound}
\mathrm{Pr}\left( \lambda_{\max}\left( \sum\limits_{i=1}^{n} \mathcal{X}_i \right) \geq \alpha \right) \leq mp  \inf\limits_{t > 0}  \exp(- t \alpha + (e^t - 1) \mu_{\max})
\end{eqnarray}
By selecting $t = \log (1 + \theta)$ and $\alpha \rightarrow (1 + \theta) \mu_{\max}$, we can establish  Eq.~\eqref{eq:Chernoff II Upper Bound}.

From  Eq.~\eqref{eq1:Chernoff I Lower Bound proof} and the inequality $\log(1 + x) \leq x$ for $x > -1$, we have 
\begin{eqnarray}\label{eq1:Chernoff II Lower Bound}
\mathrm{Pr}\left( \lambda_{\min}\left( \sum\limits_{i=1}^{n} \mathcal{X}_i \right) \leq \alpha \right) \leq mp  \inf\limits_{t > 0} \exp(- t \alpha - (e^t - 1) \mu_{\min})
\end{eqnarray}
By selecting $t =- \log (1 - \theta)$ and $\alpha \rightarrow (1 - \theta) \mu_{\min}$, we can establish  Eq.~\eqref{eq:Chernoff II Lower Bound}. Therefore, this theorem is proved. 
$\hfill \Box$

\subsection{T-product tensor Chernoff Inequalities for Eigentuple}

In this section, we wll present T-product tensor Chernoff inequalities about the maximum of eigentuple. 

\TensorChernoffBoundIEigentuple*

\textbf{Proof:}
From Lemma~\ref{lma:Chernoff MGF}, we have 
\begin{eqnarray}
\mathcal{I} + f(t) \mathbb{E} \mathcal{X}_i \succeq \mathbb{E} e^{t \mathcal{X}_i},
\end{eqnarray}
where $f(t) \define e^t -1$ for $t > 0$.
By applying Corollary~\ref{cor:3_9 eigentuple}, we obtain
\begin{eqnarray}\label{eq1:Chernoff I Upper Bound proof eigentuple}
\mathrm{Pr} \left( \mathbf{d}_{\max}\left( \sum\limits_{i=1}^{n} \mathcal{X}_i \right) \geq \mathbf{b} \right) &\leq& mp  \inf\limits_{t > 0} \min\limits_{1 \leq j \leq p} \left\{
 \frac{
\exp \left( n \log \lambda_{\max}\left( \frac{1}{n} \sum\limits_{i=1}^{n}\mathbb{E}e^{t \mathcal{X}_j } \right) \right) 
      }{ \left( e_{\bigodot}^{t \mathbf{b} }\right)_j }
 \right\} \nonumber \\
&\leq &  mp  \inf\limits_{t > 0} 
\exp \left( - t b_{\tilde{j}}+ n \log \lambda_{\max}\left( \frac{1}{n} \sum\limits_{i=1}^{n}\mathbb{E}e^{t \mathcal{X}_i } \right) \right) \nonumber \\
&\leq &  mp  \inf\limits_{t > 0} 
\exp \left( - t b_{\tilde{j}}+ n \log \lambda_{\max}\left(
\mathcal{I}+f(t) \frac{1}{n} \sum\limits_{i=1}^{n}\mathbb{E}\mathcal{X}_i \right) \right)  \nonumber \\
& = &  mp  \inf\limits_{t > 0} 
\exp \left( - t b_{\tilde{j}}+ n \log \left(1 + f(t) \overline{\mu}_{\max} \right) \right),
\end{eqnarray}
The last equality follows from the definition of $\overline{\mu}_{\max}$ and spectral mapping theorem. When the value $t$ at the right-hand side of  Eq.~\eqref{eq1:Chernoff I Upper Bound proof eigentuple} is 
\begin{eqnarray}\label{eq2:Chernoff I Upper Bound proof eigentuple}
t = \log \frac{b_{\tilde{j}} }{n - b_{\tilde{j}}} - \log \frac{\overline{\mu}_{\max}}{1 - \overline{\mu}_{\max}},
\end{eqnarray}
we can achieve the tightest upper bound at Eq.~\eqref{eq1:Chernoff I Upper Bound proof eigentuple}. By substituting the value $t$ in Eq.~\eqref{eq2:Chernoff I Upper Bound proof eigentuple} into Eq.~\eqref{eq1:Chernoff I Upper Bound proof eigentuple}, Eq.~\eqref{eq:Chernoff I Upper Bound eigentuple} is proved. The next goal is to prove Eq.~\eqref{eq:Chernoff I Lower Bound eigentuple}.

If we apply Lemma~\ref{lma:Chernoff MGF} to the sequence $\{- \mathcal{X}_i \}$, we have 
\begin{eqnarray}
\mathcal{I} - g(t) \mathbb{E} \mathcal{X}_i \succeq \mathbb{E} e^{t (-\mathcal{X}_i)},
\end{eqnarray}
where $g(t) \define 1 - e^t$ for $t > 0$.
By applying Corollary~\ref{cor:3_9 eigentuple} again, we obtain
\begin{eqnarray}\label{eq1:Chernoff I Lower Bound proof eigentuple}
\mathrm{Pr} \left( \lambda_{\min}\left( \sum\limits_{i=1}^{n} \mathcal{X}_i \right) \leq \mathbf{b} \right) &=& \mathrm{Pr} \left( \lambda_{\max}\left( \sum\limits_{i=1}^{n} \left( - \mathcal{X}_i \right)\right) \geq \mathbf{b} \right) 
\nonumber \\
&\leq&  mp  \inf\limits_{t > 0} \min\limits_{1 \le j \le p} \left\{
 \frac{
\exp \left( n \log \lambda_{\max}\left( \frac{1}{n} \sum\limits_{i=1}^{n}\mathbb{E}e^{ - t \mathcal{X}_i } \right) \right) 
      }{ \left( e_{\bigodot}^{t \mathbf{b} }\right)_j }
 \right\}
\nonumber \\
&\leq& mp \exp\left( -t b_{\tilde{j}} + n \log\lambda_{\max}\left( \frac{1}{n} \sum\limits_{i=1}^{n} \left( \mathcal{I} - g(t) \mathbb{E} \mathcal{X}_i \right) \right) \right) \nonumber \\
&=_1 &  mp \exp\left( -t b_{\tilde{j}}  + n \log \left( 1 -f(t)  \lambda_{\min}\left(    \frac{1}{n} \sum\limits_{i=1}^{n} \mathbb{E} \mathcal{X}_i \right) \right) \right) \nonumber \\ 
&=& mp \exp\left(-t b_{\tilde{j}} + n \log \left(1 - g(t) \overline{\mu}_{\min} \right) \right),
\end{eqnarray}
where we apply the relation $\lambda_{\min}( -\frac{1}{n} \sum\limits_{i=1}^{n} \mathbb{E} \mathcal{X}_i  ) = -\lambda_{\max} (\frac{1}{n} \sum\limits_{i=1}^{n} \mathbb{E} \mathcal{X}_i )$ at the equality $=_1$. When the value $t$ at the right-hand side of  Eq.~\eqref{eq1:Chernoff I Lower Bound proof eigentuple} is 
\begin{eqnarray}\label{eq2:Chernoff I Lower Bound proof eigentuple}
t = \log \frac{\overline{\mu}_{\min}}{1 - \overline{\mu}_{\min}} -  \log \frac{b_{\tilde{j}} }{n - b_{\tilde{j}}} ,
\end{eqnarray}
we can achieve the tightest upper bound at Eq.~\eqref{eq1:Chernoff I Lower Bound proof eigentuple}. By substituting the value $t$ in Eq.~\eqref{eq2:Chernoff I Lower Bound proof eigentuple} into Eq.~\eqref{eq1:Chernoff I Lower Bound proof eigentuple}, therefore,  Eq.~\eqref{eq:Chernoff I Lower Bound eigentuple} is proved also.
$\hfill \Box$

The tensor Chernoff bounds discussed at Theorem~\ref{thm:TensorChernoffBoundI for Eigentuple} is not related to $\mu_{\max}$ and $\mu_{\min}$ directly. Following theorem is another version of tensor Chernoff bounds to associate the probability range in terms of $\mu_{\max}$ and $\mu_{\min}$ directly and these formats of tensor Chernoff bounds are easier to be applied. 

\TensorChernoffBoundIIEigentuple*

%
\textbf{Proof:}
Without loss of generality, we can assume $T=1$ in our proof. From  Eq.~\eqref{eq1:Chernoff I Upper Bound proof eigentuple} and the inequality $\log(1 + x) \leq x$ for $x > -1$, we have 
\begin{eqnarray}\label{eq1:Chernoff II Upper Bound eigentuple}
\mathrm{Pr}\left( \mathbf{d}_{\max}\left( \sum\limits_{i=1}^{n} \mathcal{X}_i \right) \geq (1+\theta) \mu_{\max} \mathbf{1} \right) \leq mp  \inf\limits_{t > 0} \exp(- t b_{\tilde{j}} + (e^t - 1)\mu_{\max})
\end{eqnarray}
By selecting $t = \log (1 + \theta)$ and $b_{\tilde{j}} \rightarrow (1 + \theta) \mu_{\max}$, we can establish  Eq.~\eqref{eq:Chernoff II Upper Bound eigentuple}.

From  Eq.~\eqref{eq1:Chernoff I Lower Bound proof eigentuple} and the inequality $\log(1 + x) \leq x$ for $x > -1$, we have 
\begin{eqnarray}\label{eq1:Chernoff II Lower Bound eigentuple eigentuple}
\mathrm{Pr}\left( \lambda_{\min}\left( \sum\limits_{i=1}^{n} \mathcal{X}_i \right) \leq (1-\theta) \mu_{\min} \mathbf{1} \right) \leq mp \inf\limits_{t > 0} \exp(- t b_{\tilde{j}} -  (e^t - 1) \mu_{\min})
\end{eqnarray}
By selecting $t =- \log (1 - \theta)$ and $b_{\tilde{j}} \rightarrow (1 - \theta) \mu_{\min}$, we can establish  Eq.~\eqref{eq:Chernoff II Lower Bound eigentuple}. Therefore, this theorem is proved. 
$\hfill \Box$

\subsection{Application of T-product Tensor Chernoff Bounds}\label{sec:Application of Tensor Chernoff Bounds}

One application of T-product tensor Chernoff bounds is to estimate the expectation of the maximum eigenvalue of independent sum of random T-product tensors. 
\begin{corollary}[Upper and Lower Bounds for the Maximum Eigenvalue]\label{cor:Bounds for the Maximum Eigenvalue}
Consider a sequence $\{ \mathcal{X}_i  \in \mathbb{C}^{m \times m \times p } \}$ of independent, random, Hermitian T-product tensors that satisfy
\begin{eqnarray}
\mathcal{X}_i \succeq \mathcal{O} \mbox{~~and~~} \lambda_{\max}(\mathcal{X}_i) \leq T
\mbox{~~ almost surely.}
\end{eqnarray}
Then, we have
\begin{eqnarray}\label{eq1:cor:Bounds for the Maximum Eigenvalue}
\mu_{\max} \leq \mathbb{E} \lambda_{\max}\left( \sum\limits_{i=1}^n \mathcal{X}_i \right) \leq C mp e^{- \mu_{\max}/T},
\end{eqnarray}
where the constant value of $C$ is about 10.28.
\end{corollary}
\textbf{Proof:}
The lower bound at Eq.~\eqref{eq1:cor:Bounds for the Maximum Eigenvalue} is true from the convexity of the function $\mathcal{A} \rightarrow \lambda_{\max}(\mathcal{A})$ and the Jensen's inequality. 

For the upper bound, we have 
\begin{eqnarray}\label{eq2:cor:Bounds for the Maximum Eigenvalue}
\mathbb{E} \lambda_{\max}\left( \sum\limits_{i=1}^n \mathcal{X}_i \right)
&=& \int_{0}^{\infty} \mathrm{Pr} \left( \lambda_{\max}\left( \sum\limits_{i=1}^n \mathcal{X}_i \right)  \geq t \right) d t  \nonumber \\
&\leq_1& \int_{0}^{\infty} mp \exp(- \delta t + (e^{\delta} - 1) \mu_{\max}/T)dt \nonumber \\
&  &~~\mbox{$\delta$ is a positive real variable to be optimized} \nonumber \\
&=&  \frac{e^{e^{\delta}}}{\delta} mp e^{- \mu_{\max}/T} \nonumber \\
&\leq&  \frac{e^{e^{\delta_{opt}}}}{\delta_{opt}} mp e^{- \mu_{\max}/T}
 = C mp e^{- \mu_{\max}/T},
\end{eqnarray}
where the inequality $\leq_1$ comes from  Eq.~\eqref{eq1:Chernoff II Upper Bound} with the scaling factor $T$. If we select $\theta$ as the solution of the following relation $e^{\delta_{opt}} = \frac{1}{\delta_{opt}}$ to minimize the right-hand side of Eq.~\eqref{eq2:cor:Bounds for the Maximum Eigenvalue},  we have the desired upper bound when $\delta_{opt} \approx = 0.56699$. This corollary is proved. 
$\hfill \Box$

\section{Bernstein Bounds for T-product Tensors}\label{sec:Bernstein Bounds for T-Symmetric Tensors}

For random variables, Bernstein inequalities give the upper tail of a sum of
independent, zero-mean random variables that are either bounded or subexponential. In this section, we wish to extend Bernstein bounds for a sum of zero-mean random T-product tensors.

\subsection{T-product Tensor Bernstein Bounds Derivation}\label{sec:Tensor Bernstein Bounds Derivation}

We will condier bounded T-product tensor Bernstein bounds first by considering the bounded Bernstein moment-generating function with the following Lemma.
\begin{lemma}\label{lam:Bounded Bernstein mgf}
Given a random Hermitian T-product tensor $\mathcal{X} \in \mathbb{C}^{m \times m \times p }$ that satisfies:
\begin{eqnarray}\label{eq1:lam:Bounded Bernstein mgf}
\mathbb{E} \mathcal{X} = 0 \mbox{~~and~~} \lambda_{\max}(\mathcal{X} ) \leq 1 
\mbox{~~almost surely.} 
\end{eqnarray}
Then, we have
\begin{eqnarray}
 e^{  (e^t - t - 1) \mathbb{E}(\mathcal{X}^2)}   \succeq  \mathbb{E}e^{t \mathcal{X}}
\end{eqnarray}
where $t > 0$.
\end{lemma}
\textbf{Proof:}
If we define a real function $g(x) \define \frac{e^{t x} - t x - 1}{x^2}$, it is easy to see that this function $g(x)$ is an increasing function for $0 < x \leq 1$. From Eq~\eqref{eq:transfer rule}, we have 
\begin{eqnarray}\label{eq2:lam:Bounded Bernstein mgf}
g(\mathcal{X}) \preceq g(1) \mathcal{I}.
\end{eqnarray}
Moreover, we also have
\begin{eqnarray}\label{eq3:lam:Bounded Bernstein mgf}
e^{t \mathcal{X}} &=& \mathcal{I} + t \mathcal{X} + g(\mathcal{X}) \star \mathcal{X}^2 \nonumber \\
&\preceq&  \mathcal{I} + t \mathcal{X} + g(1) \mathcal{X}^2,
\end{eqnarray}
where the $\preceq$ comes from Eq.~\eqref{eq2:lam:Bounded Bernstein mgf}. By taking the expectation for both sides of Eq.~\eqref{eq3:lam:Bounded Bernstein mgf}, we then obtain
\begin{eqnarray}
\mathbb{E}e^{t \mathcal{X}} &\preceq& \mathcal{I} + g(1) \mathbb{E}\left( \mathcal{X}^2\right) \preceq e^{g(1)\mathbb{E}\left( \mathcal{X}^2\right) } \nonumber \\
&=&   e^{ (e^t - t - 1) \mathbb{E}(\mathcal{X}^2)}.
\end{eqnarray}
This lemma is established. 
$\hfill \Box$

We are ready to present Bernstein inequalities for random T-product tensors with bounded $\lambda_{\max}$. 

\BoundedTensorBernstein*

%
\textbf{Proof:}
Without loss of generality, we can assume that $T=1$ since the summands are 1-homogeneous and the variance is 2-homogeneous. From Lemma~\ref{lam:Bounded Bernstein mgf}, we have 
\begin{eqnarray}\label{eq5:thm:Bounded Tensor Bernstein}
\mathbb{E}e^{t \mathcal{X}_i} \preceq e^{(e^t - t- 1) \mathbb{E}(\mathcal{X}_i^2)} \mbox{~~ for $t > 0$.}
\end{eqnarray}
By applying Corollary~\ref{cor:3_7}, we then have 
\begin{eqnarray}\label{eq6:thm:Bounded Tensor Bernstein}
\mathrm{Pr} \left( \lambda_{\max}\left( \sum\limits_{i=1}^{n} \mathcal{X}_i \right)\geq \theta \right) &\leq &mp  \exp\left(- t \theta + (e^t - t - 1) \lambda_{\max}\left( \sum\limits_{i=1}^n \mathbb{E} \left( \mathcal{X}^2_i \right)\right)\right) \nonumber \\
&=&mp \exp \left( - t \theta + \sigma^2(e^t - t - 1)   \right).
\end{eqnarray}
The right-hand side of Eq.~\eqref{eq6:thm:Bounded Tensor Bernstein} can be minimized by setting $t = \log (1 + \theta/\sigma^2)$. Substitute such $t$ and simplify the right-hand side of Eq.~\eqref{eq6:thm:Bounded Tensor Bernstein}, we obtain  Eq.~\eqref{eq2:thm:Bounded Tensor Bernstein}.

For $\theta \leq \sigma^2/T$, we have 
\begin{eqnarray}
\frac{1}{\sigma^2 + T\theta/3} \geq \frac{1}{\sigma^2 + T (\sigma^2/T) /3} = \frac{3}{4\sigma^2}, 
\end{eqnarray}
then, we obtain  Eq.~\eqref{eq3:thm:Bounded Tensor Bernstein}. Correspondingly, for $\theta \geq \sigma^2/T$, we have 
\begin{eqnarray}
\frac{\theta}{\sigma^2 + T\theta/3} \geq \frac{\sigma^2/T}{\sigma^2 + T (\sigma^2/T)/3} = \frac{3}{4 T}, 
\end{eqnarray}
and, we obtain  Eq.~\eqref{eq4:thm:Bounded Tensor Bernstein} also. 
$\hfill \Box$

The following theorem~\ref{thm:Subexponential Tensor Bernstein} is the extension of the theorem~\ref{thm:Bounded Tensor Bernstein} by allowing the moments of the random T-product tensors to grow at a controlled rate. We have to prepare subexponential Bernstein moment-generating function Lemma first for later proof of Theorem~\ref{thm:Subexponential Tensor Bernstein}

\begin{lemma}\label{lem:Subexponential Bernstein mgf}
Suppose that $\mathcal{X}$ is a random Hermitian T-product tensor that satisfies
\begin{eqnarray}\label{eq1:lem:Subexponential Bernstein mgf}
\mathbb{E} \mathcal{X} = 0 \mbox{~~and~~} \mathbb{E}(\mathcal{X}^p)  \preceq \frac{p! \mathcal{A}^2}{2} 
~~ \mbox{for $p=2,3,4,\cdots$.} 
\end{eqnarray}
Then, we have
\begin{eqnarray}
 \exp\left(\frac{t^2 \mathcal{A}^2}{2(1 - t)}\right)   \succeq  \mathbb{E}e^{t \mathcal{X}},
\end{eqnarray}
where $0 < t < 1$.
\end{lemma}
\textbf{Proof:}
From Taylor series of the tensor exponential expansion, we have 
\begin{eqnarray}
\mathbb{E}e^{t \mathcal{X}} &=& \mathcal{I} + t \mathbb{E}\mathcal{X} + \sum\limits_{p=2}^{\infty} \frac{t^p \mathbb{E}(\mathcal{X}^p)}{p!} \preceq \mathcal{I}+ \sum\limits_{p=2}^{\infty} \frac{t^p \mathcal{A}^2}{2} \nonumber \\
&=& \mathcal{I} + \frac{t^2 \mathcal{A}^2}{2 (1-t)} \preceq \exp\left(\frac{t^2 \mathcal{A}^2}{2(1 - t)} \right),
\end{eqnarray}
therefore, this Lemma is proved. 
$\hfill \Box$

\SubexponentialTensorBernstein*

%
%
\textbf{Proof:}
Without loss of generality, we can assume that $T=1$. From Lemma~\ref{lem:Subexponential Bernstein mgf}, we have 
\begin{eqnarray}
\mathbb{E}\exp \left(t \mathcal{X}_i \right) \preceq \exp\left(\frac{t^2 \mathcal{A}^2_i}{2(1-t)} \right), 
\end{eqnarray}
where $0 < t < 1$.

By applying Corollary~\ref{cor:3_7}, we then have 
\begin{eqnarray}\label{eq6:thm:Subexponential Tensor Bernstein}
\mathrm{Pr} \left( \lambda_{\max}\left( \sum\limits_{i=1}^{n} \mathcal{X}_i \right)\geq \theta \right) &\leq & mp \exp\left(- t \theta + \frac{t^2}{2(1-t)} \lambda_{\max}\left( \sum\limits_{i=1}^n \mathcal{A}_i^2  \right)\right) \nonumber \\
&=& mp \exp \left( - t \theta + \frac{\sigma^2 t^2}{2(1-t)}  \right).
\end{eqnarray}
The right-hand side of Eq.~\eqref{eq6:thm:Subexponential Tensor Bernstein} can be minimized by setting $t = \frac{\theta}{\theta + \sigma^2}$. Substitute such $t$ and simplify the right-hand side of Eq.~\eqref{eq6:thm:Subexponential Tensor Bernstein}, we obtain  Eq.~\eqref{eq2:thm:Subexponential Tensor Bernstein}.

For $\theta \leq \sigma^2/T$, we have 
\begin{eqnarray}
\frac{1}{\sigma^2 + T\theta} \geq \frac{1}{\sigma^2 + T (\sigma^2/T )} = \frac{1}{2\sigma^2}, 
\end{eqnarray}
then, we obtain  Eq.~\eqref{eq3:thm:Subexponential Tensor Bernstein}. Similarly, for $\theta \geq \sigma^2/T$, we have 
\begin{eqnarray}
\frac{\theta}{\sigma^2 + T\theta} \geq \frac{\sigma^2/T}{\sigma^2 + T( \sigma^2/T) } = \frac{1}{2T}, 
\end{eqnarray},
therefore, we also obtain  Eq.~\eqref{eq4:thm:Subexponential Tensor Bernstein}. 
$\hfill \Box$

\subsection{T-product Tensor Bernstein Bounds for Eigentuple}\label{sec:Tensor Bernstein Bounds Derivation eigentuple}

In this section, we will extend T-product tensor bernstein bounds from the maximum eigenvalue discussed at previous section to the maximum eigentuple.

\BoundedTensorBernsteinEigentuple*

\textbf{Proof:}
Without loss of generality, we can assume that $T=1$ since the summands are 1-homogeneous and the variance is 2-homogeneous. From Lemma~\ref{lam:Bounded Bernstein mgf}, we have 
\begin{eqnarray}\label{eq5:thm:Bounded Tensor Bernstein  eigentuple}
\mathbb{E}e^{t \mathcal{X}_i} \preceq e^{(e^t - t- 1) \mathbb{E}(\mathcal{X}_i^2)} \mbox{~~ for $t > 0$.}
\end{eqnarray}
By applying Corollary~\ref{cor:3_7 eigentuple}, we then have 
\begin{eqnarray}\label{eq6:thm:Bounded Tensor Bernstein eigentuple}
\mathrm{Pr} \left( \mathbf{d}_{\max}\left( \sum\limits_{i=1}^{n} \mathcal{X}_i \right)\geq \mathbf{b} \right) &\leq &  mp \inf\limits_{t > 0} \min\limits_{1 \leq j \leq p} \left\{ \frac{ \exp  \left(  f(t) \lambda_{\max} \left( \sum\limits_{i=1}^n \mathcal{A}_i \right)    \right)    }{ \left( e_{\bigodot}^{t \mathbf{b} }\right)_j } \right\} \nonumber \\
& \leq & mp  \exp\left(- t b_{\tilde{j}} + (e^t - t - 1) \lambda_{\max}\left( \sum\limits_{i=1}^n \mathbb{E} \left( \mathcal{X}^2_i \right)\right)\right) \nonumber \\
&=&mp \exp \left( - t b_{\tilde{j}}  + \sigma^2(e^t - t - 1)   \right).
\end{eqnarray}
The right-hand side of Eq.~\eqref{eq6:thm:Bounded Tensor Bernstein eigentuple} can be minimized by setting $t = \log (1 + b_{\tilde{j}} /\sigma^2)$. Substitute such $t$ and simplify the right-hand side of Eq.~\eqref{eq6:thm:Bounded Tensor Bernstein eigentuple}, we obtain  Eq.~\eqref{eq2:thm:Bounded Tensor Bernstein eigentuple}.

For $b_{\tilde{j}} \leq \sigma^2/T$, we have 
\begin{eqnarray}
\frac{1}{\sigma^2 + Tb_{\tilde{j}}/3} \geq \frac{1}{\sigma^2 + T (\sigma^2/T) /3} = \frac{3}{4\sigma^2}, 
\end{eqnarray}
then, we obtain  Eq.~\eqref{eq3:thm:Bounded Tensor Bernstein eigentuple}. Correspondingly, for $b_{\tilde{j}} \geq \sigma^2/T$, we have 
\begin{eqnarray}
\frac{\theta}{\sigma^2 + Tb_{\tilde{j}}/3} \geq \frac{\sigma^2/T}{\sigma^2 + T (\sigma^2/T)/3} = \frac{3}{4 T}, 
\end{eqnarray}
and, we obtain  Eq.~\eqref{eq4:thm:Bounded Tensor Bernstein eigentuple} also. 
$\hfill \Box$

Below theorem is another variation of T-product tensor Bernstein bounds by subexponential 
constraints of $\mathbb{E} (\mathcal{X}^p_i)$. 

\SubexponentialTensorBernsteinEigentuple*

\textbf{Proof:}
Without loss of generality, we can assume that $T=1$. From Lemma~\ref{lem:Subexponential Bernstein mgf}, we have 
\begin{eqnarray}
\mathbb{E}\exp \left(t \mathcal{X}_i \right) \preceq \exp\left(\frac{t^2 \mathcal{A}^2_i}{2(1-t)} \right), 
\end{eqnarray}
where $0 < t < 1$.

By applying Corollary~\ref{cor:3_7 eigentuple}, we then have 
\begin{eqnarray}\label{eq6:thm:Subexponential Tensor Bernstein eigentuple}
\mathrm{Pr} \left( \lambda_{\max}\left( \sum\limits_{i=1}^{n} \mathcal{X}_i \right)\geq \theta \right) & \leq & mp \inf\limits_{t > 0} \min\limits_{i} \left\{ \frac{ \exp  \left(  f(t) \lambda_{\max} \left( \sum\limits_{i=1}^n \mathcal{A}_i \right)    \right)    }{ \left( e_{\bigodot}^{t \mathbf{b} }\right)_i } \right\} \nonumber \\
& \leq & 
mp \exp\left(- t b_{\tilde{j}} + \frac{t^2}{2(1-t)} \lambda_{\max}\left( \sum\limits_{i=1}^n \mathcal{A}_i^2  \right)\right) \nonumber \\
&=& mp \exp \left( - t b_{\tilde{j}}+ \frac{\sigma^2 t^2}{2(1-t)}  \right).
\end{eqnarray}
The right-hand side of Eq.~\eqref{eq6:thm:Subexponential Tensor Bernstein eigentuple} can be minimized by setting $t = \frac{b_{\tilde{j}}}{b_{\tilde{j}}+ \sigma^2}$. Substitute such $t$ and simplify the right-hand side of Eq.~\eqref{eq6:thm:Subexponential Tensor Bernstein eigentuple}, we obtain  Eq.~\eqref{eq2:thm:Subexponential Tensor Bernstein eigentuple}.

For $b_{\tilde{j}} \leq \sigma^2/T$, we have 
\begin{eqnarray}
\frac{1}{\sigma^2 + Tb_{\tilde{j}}} \geq \frac{1}{\sigma^2 + T (\sigma^2/T )} = \frac{1}{2\sigma^2}, 
\end{eqnarray}
then, we obtain  Eq.~\eqref{eq3:thm:Subexponential Tensor Bernstein eigentuple}. Similarly, for $b_{\tilde{j}} \geq \sigma^2/T$, we have 
\begin{eqnarray}
\frac{\theta}{\sigma^2 + Tb_{\tilde{j}}} \geq \frac{\sigma^2/T}{\sigma^2 + T( \sigma^2/T) } = \frac{1}{2T}, 
\end{eqnarray},
therefore, we also obtain  Eq.~\eqref{eq4:thm:Subexponential Tensor Bernstein eigentuple}. 
$\hfill \Box$

\subsection{Application of Tensor Bernstein Bounds}\label{sec:Application of Tensor Bernstein Bounds}

The tensor Bernstein bounds can also be extended to rectangular tensors by dilation. Consider a sequence of tensors $\{\mathcal{Y}_i\} \in \mathbb{C}^{m \times n \times p}$ satisfy following:
\begin{eqnarray}
\mathbb{E} \mathcal{Y}_i = \mathcal{O} \mbox{~~and~~} \left\Vert \mathcal{Y}_i \right\Vert \leq T \mbox{~~almost surely.}
\end{eqnarray}
If the variance $\sigma^2$ is expressed as:
\begin{eqnarray}
\sigma^2 \define  \max \Bigg\{ \left\Vert \sum\limits_{i=1}^n \mathcal{Y}_i \star \mathcal{Y}^H_i \right\Vert,  \left\Vert \sum\limits_{i=1}^n \mathcal{Y}^H_i \star \mathcal{Y}_i \right\Vert\Bigg\},
\end{eqnarray}
we have 
\begin{eqnarray}
\mathrm{Pr} \left( \left\Vert \sum\limits_{i=1}^{n} \mathcal{Y}_i \right\Vert \geq \theta \right) \leq (m+n) p \exp \left( \frac{-\theta^2/2}{\sigma^2 + T\theta/3}\right)
\end{eqnarray}
from Theorem~\ref{thm:Bounded Tensor Bernstein}. 

Another application of tensor Bernstein bounds is to get upper and lower Bounds for the maximum eigenvalue with subexponential tensors. This application can relax the corollary~\ref{cor:Bounds for the Maximum Eigenvalue} conditions by allowing the moments of the random tensors to grow at a controlled rate. 
\begin{corollary}[Upper and Lower Bounds for the Maximum Eigenvalue for Subexponential Tensors]\label{cor:Bounds for the Maximum Eigenvalue Subexponential Tensors}
Consider a sequence $\{ \mathcal{X}_i  \in \mathbb{C}^{m \times m \times p} \}$ of independent, random, Hermitian T-product tensors that satisfy
\begin{eqnarray}
\mathcal{X}_i \succeq \mathcal{O} \mbox{~~and~~} \mathbb{E} (\mathcal{X}^p_i) \preceq \frac{p! T^{p-2} }{2} \mathcal{A}_i^2,
\end{eqnarray}
and $\sigma^2 \define \left\Vert \sum\limits_{i=1}^n \mathcal{A}^2_i \right\Vert$. Then, we have 
\begin{eqnarray}\label{eq1:cor:Bounds for the Maximum Eigenvalue Subexponential Tensors}
\mu_{\max} \leq \mathbb{E} \lambda_{\max}\left( \sum\limits_{i=1}^n \mathcal{X}_i \right) \leq
2mp\left(\sigma \mathfrak{I}\left(\frac{\sigma}{2T} \right) + 2T e^{-\frac{\sigma^2}{4T^2}}\right),
\end{eqnarray}
where $\mathfrak{I}(\frac{\sigma}{2T}) \define \int_0^{\frac{\sigma}{2T}} e^{-s^2} ds$.
\end{corollary}
\textbf{Proof:}
The lower bound at Eq.~\eqref{eq1:cor:Bounds for the Maximum Eigenvalue Subexponential Tensors} is true from the convexity of the function $\mathcal{A} \rightarrow \lambda_{\max}(\mathcal{A})$ and the Jensen's inequality. 

For the upper bound, we have 
\begin{eqnarray}\label{eq2:cor:Bounds for the Maximum Eigenvalue Subexponential Tensors}
\mathbb{E} \lambda_{\max}\left( \sum\limits_{i=1}^n \mathcal{X}_i \right)
&=& \int_{0}^{\infty} \mathrm{Pr} \left( \lambda_{\max}\left( \sum\limits_{i=1}^n \mathcal{X}_i \right)  \geq t \right) d t  \nonumber \\
&\leq_1& mp \int_{0}^{\frac{\sigma^2}{T}}\exp \left(- \frac{t^2}{4 \sigma^2} \right)   dt +  mp  \int_{\frac{\sigma^2}{T}}^{\infty}\exp \left(- \frac{t}{4 T} \right)   dt \nonumber \\
&=&  2mp\left(\sigma \mathfrak{G}\left(\frac{\sigma}{2T} \right) + 2T e^{-\frac{\sigma^2}{4T^2}}\right),
\end{eqnarray}
where the inequality $\leq_1$ comes from the Eqs.~\eqref{eq3:thm:Subexponential Tensor Bernstein} and~\eqref{eq4:thm:Subexponential Tensor Bernstein}. This corollary is proved by introducing Gaussian integral function $\mathfrak{G}(x) \define \int_0^{x} e^{-s^2} ds$.  
$\hfill \Box$

\section{T-product Tensor Martingales Inequalities}\label{sec:T-Symmetric Tensor Martingales Inequalities}

In this section, we introduce concepts about T-product tensor martingales in Section~\ref{sec:Tensor Martingales}, and extend Hoeffding, Azuma, and McDiarmid inequalities to random T-product tensors context in Section~\ref{sec:Tensor Martingale Deviation Bounds}. These bounds are extended to the eigentuple version in Section~\ref{sec:Tensor Martingale Deviation Bounds Eigentuple}.  

\subsection{T-product Tensor Martingales}\label{sec:Tensor Martingales}

Several basic definitions about T-product tensor martingales will be provided here for later T-product tensor martingales related bounds. Let $(\Omega, \mathfrak{F}, \mathbb{P})$ be a master probability space. Consider a filtration $\{ \mathfrak{F}_i \}$ contained in the master sigma algebra as:
\begin{eqnarray}
\mathfrak{F}_0 \subset \mathfrak{F}_1 \subset \mathfrak{F}_2 \subset \cdots 
\subset \mathfrak{F}_{\infty} \subset \mathfrak{F}.
\end{eqnarray}
Given such a filtration, we define the conditional expectation $\mathbb{E}_i[ \cdot ] \define \mathbb{E}_i[ \cdot | \mathfrak{F}_i]$. A sequence $\{ \mathcal{Y}_i \}$ of random tensors is called \emph{adapted} to the filtration when each tensor $\mathcal{Y}_i$ is measurable with respect to $\mathfrak{F}_i$. We can think that an adapted sequence is one where the present depends only on the past. 

An adapted sequence $\{ \mathcal{X}_i \}$ of Hermitian T-product tensors is named as a \emph{tensor martingale} when 
\begin{eqnarray}
\mathbb{E}_{i-1} \mathcal{X}_i = \mathcal{X}_{i-1} \mbox{~~~and~~~} \mathbb{E} \left\Vert \mathcal{X}_i \right\Vert < \infty,
\end{eqnarray}
where $i = 1, 2, 3, \cdots$. We obtain a scalar martingale if we track any fixed entry of a tensor martingale $\{ \mathcal{X}_i \}$. Given a tensor martingale $\{ \mathcal{X}_i \}$, we can construct the following new sequence of tensors 
\begin{eqnarray}
\mathcal{Y}_i \define \mathcal{X}_i  - \mathcal{X}_{i-1}  \mbox{~~for $i=1, 2, 3,\cdots$}
\end{eqnarray}
We then have $\mathbb{E}_{i-1}\mathcal{Y}_i = \mathcal{O}$.

\subsection{Tensor Martingale Deviation Bounds for Eigenvalues}\label{sec:Tensor Martingale Deviation Bounds}

Two Lemmas should be presented first before presenting tensor martingale deviation bounds and their proofs.

\begin{lemma}[Tensor Symmetrization]\label{lma:tensor symmetriation}
Let $\mathcal{A} \in \mathbb{C}^{m \times m \times p}$ be a fixed Hermitian T-product tensor, and let $\mathcal{X}$ be a random Hermitian T-product tensor with $\mathbb{E}\mathcal{X} = \mathcal{O}$. Then
\begin{eqnarray}
\mathbb{E} \mathrm{Tr}e^{\mathcal{A} + \mathcal{X}} \leq \mathbb{E} \mathrm{Tr}e^{\mathcal{A} + 2\beta \mathcal{X}},
\end{eqnarray}
where $\beta$ is a Rademacher random variable.
\end{lemma}
\textbf{Proof:}
Build an independent copy random tensor $\mathcal{Y}$ from $\mathcal{X}$, and let $\mathbb{E}_{\mathcal{Y}}$ denote the expectation with respect to the new random tensor $\mathcal{Y}$. Then, we have
\begin{eqnarray}
\mathbb{E} \mathrm{Tr}e^{\mathcal{A} + \mathcal{X}} =
\mathbb{E} \mathrm{Tr}e^{\mathcal{A} + \mathcal{X} - \mathbb{E}_{\mathcal{Y}} \mathcal{Y}}  \leq \mathbb{E} \mathrm{Tr}e^{\mathcal{A} + \mathcal{X} - \mathcal{Y}}
= \mathbb{E} \mathrm{Tr}e^{\mathcal{A} + \beta( \mathcal{X} - \mathcal{Y})},
\end{eqnarray}
where the first equality uses $\mathbb{E}_{\mathcal{Y}} \mathcal{Y} = \mathcal{O}$; the inequality uses the convexity of the trace exponential with Jensen's inequality; finally, the last equality comes from that the random tensor $\mathcal{X} - \mathcal{Y}$ is a symmetric random tensor and Rademacher is also a symmetric random variable.

This Lemma is established by the following:
\begin{eqnarray}
\mathbb{E} \mathrm{Tr}e^{\mathcal{A} + \mathcal{X}} &\leq& \mathbb{E} \mathrm{Tr} \left( e^{\mathcal{A}/2 + \beta \mathcal{X}}  e^{\mathcal{A}/2 - \beta \mathcal{Y}}\right) \nonumber \\
& \leq & \left( \mathbb{E} \mathrm{Tr} e^{\mathcal{A} + 2 \beta \mathcal{X}}\right)^{1/2}  \left( \mathbb{E} \mathrm{Tr} e^{\mathcal{A} - 2 \beta \mathcal{Y}}\right)^{1/2}
= \mathbb{E} \mathrm{Tr}e^{\mathcal{A} + 2 \beta \mathcal{X}},
\end{eqnarray}
where the first inequality comes from T-product tensor Golden-Thompson inequality by Theorem~\ref{thm:GT Inequality for T-product Tensors}, the second inequality comes from the Cauchy-Schwarz inequality, and the last identity follows from that the two factors are identically distributed. 
$\hfill \Box$

Following lemma is introduced to provide the tensor cumulant-generating function of a symmetrized random tensor. 
\begin{lemma}[Cumulant-Generating Function of Symetrized Random T-product Tensor]\label{lma:Azuma CGF}
Given that $\mathcal{X} \in \mathbb{C}^{m \times m \times p}$ is a random Hermitian T-product tensor and $\mathcal{A} \in \mathbb{C}^{m \times m \times p}$ is a fixed Hermitian T-product tensor that satisfies $\mathcal{X}^2 \preceq \mathcal{A}^2$. Then, we have
\begin{eqnarray}
\log \mathbb{E} \left [ e^{2 \beta t  \mathcal{X}} | \mathcal{X} \right] \preceq 2 t^2 \mathcal{A}^2, 
\end{eqnarray}
where $\beta$ is a Rademacher random variable.
\end{lemma}
\textbf{Proof:}
From Lemma~\ref{lma:mgf of Rademacher and normal rvs.}, we have 
\begin{eqnarray}
\mathbb{E} \left [ e^{2 \beta t \mathcal{X}} | \mathcal{X} \right] \preceq  e^{2 t ^2 \mathcal{X}^2}.
\end{eqnarray}
And, from the monotone property of logarithm, we also have 
\begin{eqnarray}
\log \mathbb{E} \left [ e^{2 t \theta \mathcal{X}} | \mathcal{X} \right] \preceq  2t^2 \mathcal{X}^2 \preceq 2 t^2 \mathcal{A}^2~~\mbox{for $t \in \mathbb{R}$}.  
\end{eqnarray}
Therefore, this Lemma is proved. 
$\hfill \Box$

In probability theory, the Azuma inequality for a scaler martingale gives normal concentration about its mean value, and the deviation is controlled by the total maximum squared of the difference sequence. Following theorem is the T-product tensor version for Azuma inequality. 

\TensorAzuma*
%
%
\textbf{Proof:}
Define the filtration $\mathfrak{F}_i \define \mathfrak{F}(\mathcal{X}_1, \cdots, \mathcal{X}_i)$ for the process $\{\mathcal{X}_i \}$. Then, we have 
\begin{eqnarray}\label{eq3:thm:Tensor Azuma}
\mathbb{E} \mathrm{Tr} \exp\left( \sum\limits_{i=1}^n t \mathcal{X}_i \right) &=& \mathbb{E}   \left( \mathbb{E}\left( \mathrm{Tr}\exp\left( \sum\limits_{i=1}^{n-1} t \mathcal{X}_i+ t \mathcal{X}_n \right) | \mathfrak{F}_{n}\right) | \mathfrak{F}_{n-1}\right) \nonumber \\
& \leq & \mathbb{E}   \left( \mathbb{E}\left( \mathrm{Tr}\exp\left( \sum\limits_{i=1}^{n-1} t \mathcal{X}_i+ 2 \beta t \mathcal{X}_n \right) | \mathfrak{F}_{n}\right) | \mathfrak{F}_{n}\right) \nonumber \\
& \leq &  \mathbb{E}   \left( \mathrm{Tr}\exp\left( \sum\limits_{i=1}^{n-1} t \mathcal{X}_i+ \log \mathbb{E} \left( e^{2 \beta t \mathcal{X}_n} | \mathfrak{F}_{n}\right)  \right) | \mathfrak{F}_{n}\right) \nonumber \\
& \leq & \mathbb{E} \mathrm{Tr} \exp \left( \sum\limits_{i=1}^{n-1} t \mathcal{X}_i +2 t^2 \mathcal{A}_n^2 \right),
\end{eqnarray}
where the first equality comes from the total expectation property of conditional expectation; the first inequality comes from Lemma~\ref{lma:tensor symmetriation}; the second inequality comes from Corollary~\ref{cor:3.3} and the relaxation for the conditioning on the inner expectation to the larger algebra $\mathfrak{F}_n$; finally, the last inequality requires Lemma~\ref{lma:Azuma CGF}.

If we continue the iteration procedure based on Eq.~\eqref{eq3:thm:Tensor Azuma}, we have
\begin{eqnarray}\label{eq4:thm:Tensor Azuma eigenvalue}
\mathbb{E} \mathrm{Tr} \exp \left( \sum\limits_{i=1}^n t \mathcal{X}_i  \right) 
\leq \mathrm{Tr} \exp \left( 2t^2 \sum\limits_{i=1}^n \mathcal{A}^2_i \right),
\end{eqnarray}
then apply Eq.~\eqref{eq4:thm:Tensor Azuma eigenvalue} into Lemma~\ref{lma: Laplace Transform Method Eigenvalue Version}, we obtain 
\begin{eqnarray}
\mathrm{Pr} \left( \lambda_{\max}\left( \sum\limits_{i=1}^{n} \mathcal{X}_i \right)\geq \theta \right) &\leq&  \inf\limits_{t > 0} \Big\{ e^{- t \theta} \mathbb{E} \mathrm{Tr}\exp \left( \sum\limits_{i=1}^n t \mathcal{X}_i \right) \Big\} \nonumber \\
& \leq &  \inf\limits_{t > 0} \Big\{ e^{- t \theta} \mathbb{E} \mathrm{Tr}\exp \left( 2 t^2 \sum\limits_{i=1}^n  \mathcal{A}^2_i \right) \Big\} \nonumber \\
& \leq & \inf\limits_{t > 0} \Big\{ e^{- t \theta} mp 
\lambda_{\max} \left(  \exp \left( 2t^2 \sum\limits_{i=1}^{n} \mathcal{A}^2_i \right)\right) \Big\} \nonumber \\
& = &\inf\limits_{t > 0} \Big\{ e^{- t \theta} mp 
\exp \left( 2t^2 \sigma^2 \right) \Big\} \nonumber \\
& \leq &  mp  e^{-\frac{\theta^2}{8 \sigma^2}},
\end{eqnarray}
where the third inequality utilizes $\lambda_{\max}$ to bound trace, the equality applies the definition of $\sigma^2$ and spectral mapping theorem, finally, we select $t = \frac{\theta}{4 \sigma^2}$ to minimize the upper bound to obtain this theorem.
$\hfill \Box$

If we add extra assumption that the summands are independent, Theorem~\ref{thm:Tensor Azuma} gives a T-product tensor extension of Hoeffding’s inequality. If we apply Theorem~\ref{thm:Tensor Azuma} to a Hermitian T-product tensor martingale, we will have following corollary.
\begin{corollary}\label{cor:7.2}
Given a Hermitian T-product tensor martingale $\{\mathcal{Y}_i: i=1,2,\cdots,n \} \in $ \\ $\mathbb{C}^{m \times m \times p}$, and let $\mathcal{X}_i$ be the difference sequence of $\{\mathcal{Y}_i \}$, i.e., $\mathcal{X}_i \define \mathcal{Y}_i - \mathcal{Y}_{i-1}$ for $i = 1,2,3, \cdots$. If the difference sequence satisfies 
\begin{eqnarray}\label{eq1:cor:7.2}
\mathbb{E}_{i-1} \mathcal{X}_i = 0 \mbox{~~and~~} \mathcal{X}^2_i \preceq  \mathcal{A}_i ~~\mbox{almost surely}, 
\end{eqnarray}
where $i = 1,2,3,\cdots$ and the total varaince $\sigma^2$ is defined as as: $\sigma^2 \define \left\Vert \sum\limits_i^n \mathcal{A}_i^2 \right\Vert$. Then, we have
\begin{eqnarray}
\mathrm{Pr} \left( \lambda_{\max}\left( \mathcal{Y}_n -  \mathbb{E} \mathcal{Y}_n \right)\geq \theta \right) \leq mp e^{-\frac{\theta^2}{8 \sigma^2}}.
\end{eqnarray}
\end{corollary}

In the scalar setting, McDiarmid inequality can be treated as a corollary of Azuma’s inequality. McDiarmid inequality states that a function of independent random variables exhibits normal concentration about its mean, and the variance depends on the function value sensitivity with respect to the input. Following theorem is the McDiarmid inequality for the T-product tensor. 

\TensorMcDiarmid*
%
\textbf{Proof:}
We define following random tensors $\mathcal{Y}_i$ for $0 \leq i \leq n$ as:
\begin{eqnarray}\label{eq3:thm:Tensor McDiarmid}
\mathcal{Y}_i \define \mathbb{E} \left( F(x_1,\cdots,x_n) | X_1, \cdots, X_i \right) = \mathbb{E}_{X_{i+1}}  \mathbb{E}_{X_{i+2}} \cdots  \mathbb{E}_{X_{n}} F(x_1,\cdots,x_n),
\end{eqnarray}
where $ \mathbb{E}_{X_{i+1}}$ is the expectation with respect to the random variable $X_{i+1}$. The constructed sequence $\mathcal{Y}_i$ forms a martingale. The associated difference sequence with respect to $\mathcal{Y}$, denoted as $\{ \mathcal{Z}_i \}$, can be stated as:
\begin{eqnarray}\label{eq4:thm:Tensor McDiarmid}
\mathcal{Z}_i  \define \mathcal{Y}_i  - \mathcal{Y}_{i-1} = 
 \mathbb{E}_{X_{i+1}}  \mathbb{E}_{X_{i+2}} \cdots  \mathbb{E}_{X_{n}}\left(
 F(x_1,\cdots,x_n) -   \mathbb{E}_{X_{i}}  F(x_1,\cdots,x_n ) \right).
\end{eqnarray}

Because $(x_1,\cdots,x_i)$ forms a filtration with respect to $i$, we have 
\begin{eqnarray}
\mathbb{E}_{X_{i-1}}\mathcal{Y}_i &=& \mathbb{E}_{X_{i-1}} \left( \mathbb{E}_{X_{i+1}}  \mathbb{E}_{X_{i+2}} \cdots  \mathbb{E}_{X_{n}} F(x_1,\cdots,x_n) | X_{i-1} \right) \nonumber \\
&=& \mathbb{E}_{X_{i-1}} \left( \mathbb{E}_{X_{i}}  \mathbb{E}_{X_{i+1}} \cdots  \mathbb{E}_{X_{n}} F(x_1,\cdots,x_n) | X_{i-1} \right) = \mathbb{E}_{X_{i-1}}\mathcal{Y}_{i-1},
\end{eqnarray}
then,
\begin{eqnarray}\label{eq5:thm:Tensor McDiarmid}
\mathbb{E}_{X_{i-1}} \mathcal{Z}_i
&=& \mathbb{E}_{X_{i-1}}\mathcal{Y}_i - \mathbb{E}_{X_{i-1}}\mathcal{Y}_{i-1} = \mathcal{O}
\end{eqnarray}

Let $X'_i$ be an independent copy of $X_i$, and construct the following two random vectors:
\begin{eqnarray}
\mathbf{x}' &=&  (X_1,\cdots,X_{i-1},X'_i,X_{i+1},\cdots,X_n), \nonumber \\
\mathbf{x}  &=&  (X_1,\cdots,X_{i-1},X_i,X_{i+1},\cdots,X_n).
\end{eqnarray} 
Note that $\mathbb{E}_{X_i} F(\mathbf{x}) =
\mathbb{E}_{X'_i} F(\mathbf{x}')$ and $F(\mathbf{x})$ does not depend on $X'_i$, we can expresse $\mathcal{Z}_i$ from Eq.~\eqref{eq4:thm:Tensor McDiarmid} as  
\begin{eqnarray}\label{eq6:thm:Tensor McDiarmid}
\mathcal{Z}_i = \mathbb{E}_{X_{i+1}}  \mathbb{E}_{X_{i+2}} \cdots  \mathbb{E}_{X_{n}}  \mathbb{E}_{X'_{i}}  \left(
 F(\mathbf{x}) -    F(\mathbf{x}') \right).
\end{eqnarray}
Since two vectors $\mathbf{x}$ and $\mathbf{x}'$ are differ only at the $i$-th position, we have
\begin{eqnarray}
\left(  F(\mathbf{x}) -    F(\mathbf{x}')  \right)^2 \preceq \mathcal{A}^2_i,
\end{eqnarray}
from requirement provided by Eq.~\eqref{eq1:thm:Tensor McDiarmid}. Then, we have following upper bound 
\begin{eqnarray}\label{eq7:thm:Tensor McDiarmid}
\mathbb{E}_{X_{i+1}}  \mathbb{E}_{X_{i+2}} \cdots  \mathbb{E}_{X_{n}}  \mathbb{E}_{X'_{i}}  \left(
 F(\mathbf{x}) -    F(\mathbf{x}') \right)^2 \preceq \mathcal{A}^2_i.
\end{eqnarray}
Therefore, from conditions provided by Eq.~\eqref{eq5:thm:Tensor McDiarmid} and Eq.~\eqref{eq7:thm:Tensor McDiarmid}, this theorem is proved by applying Corollary~\ref{cor:7.2} to the martingale $\{\mathcal{Y}_i\}$.
$\hfill \Box$

\subsection{Tensor Martingale Deviation Bounds for Eigentuple}\label{sec:Tensor Martingale Deviation Bounds Eigentuple}

In this section, we will extend results about martingale deviation bounds for eigenvalues from Section~\ref{sec:Tensor Martingale Deviation Bounds} to martingale deviation bounds for eigentuple.

\TensorAzumaEigentuple*

%
%
\textbf{Proof:}
From Eq.~\eqref{eq3:thm:Tensor Azuma}, we have 
\begin{eqnarray}\label{eq3:thm:Tensor Azuma eigentuple}
\mathbb{E} \mathrm{Tr} \exp\left( \sum\limits_{i=1}^n t \mathcal{X}_i \right) 
& \leq & \mathbb{E} \mathrm{Tr} \exp \left( \sum\limits_{i=1}^{n-1} t \mathcal{X}_i +2 t^2 \mathcal{A}_n^2 \right).
\end{eqnarray}
If we continue the iteration procedure based on Eq.~\eqref{eq3:thm:Tensor Azuma eigentuple}, we have
\begin{eqnarray}\label{eq4:thm:Tensor Azuma eigentuple}
\mathbb{E} \mathrm{Tr} \exp \left( \sum\limits_{i=1}^n t \mathcal{X}_i  \right) 
\leq \mathrm{Tr} \exp \left( 2t^2 \sum\limits_{i=1}^n \mathcal{A}^2_i \right),
\end{eqnarray}
then apply Eq.~\eqref{eq4:thm:Tensor Azuma eigentuple} into Lemma~\ref{lma: Laplace Transform Method Eigentuple Version}, we obtain 
\begin{eqnarray}
\mathrm{Pr} \left( \mathbf{d}_{\max}\left( \sum\limits_{i=1}^{n} \mathcal{X}_i \right)\geq \mathbf{b} \right) &\leq&  
\inf_{t > 0}
 \min\limits_{1 \le j \le p} \left\{ \frac{  \mathbb{E} \left( \mathrm{Tr}\left(\exp\left( \sum\limits_{i=1}^n t \mathcal{X}_i \right) \right)\right)   }{ \left(e_{\bigodot}^{ t \mathbf{b}}  \right)_j  }\right\} \nonumber \\
& \leq &
\inf\limits_{t > 0} \Big\{ e^{- t b_{\tilde{j}}} \mathbb{E} \mathrm{Tr}\exp \left( \sum\limits_{i=1}^n t \mathcal{X}_i \right) \Big\} \nonumber \\
& \leq &  \inf\limits_{t > 0} \Big\{ e^{- t b_{\tilde{j}}} \mathbb{E} \mathrm{Tr}\exp \left( 2 t^2 \sum\limits_{i=1}^n  \mathcal{A}^2_i \right) \Big\} \nonumber \\
& \leq & \inf\limits_{t > 0} \Big\{ e^{- t b_{\tilde{j}}} mp 
\lambda_{\max} \left(  \exp \left( 2t^2 \sum\limits_{i=1}^{n} \mathcal{A}^2_i \right)\right) \Big\} \nonumber \\
& = &\inf\limits_{t > 0} \Big\{ e^{- t b_{\tilde{j}}} mp 
\exp \left( 2t^2 \sigma^2 \right) \Big\} \nonumber \\
& \leq &  mp  e^{-\frac{b_{\tilde{j}}^2}{8 \sigma^2}},
\end{eqnarray}
where the third inequality utilizes $\lambda_{\max}$ to bound trace, the equality applies the definition of $\sigma^2$ and spectral mapping theorem, finally, we select $t = \frac{b_{\tilde{j}}}{4 \sigma^2}$ to minimize the upper bound to obtain this theorem.
$\hfill \Box$

If we add an extra assumption that the summands are independent, Theorem~\ref{thm:Tensor Azuma eigentuple} gives a T-product tensor extension of Hoeffding’s inequality. If we apply Theorem~\ref{thm:Tensor Azuma eigentuple} to a Hermitian T-product tensor martingale, we will have the following corollary.
\begin{corollary}\label{cor:7.2 eigentuple}
Given a Hermitian T-product tensor martingale $\{\mathcal{Y}_i: i=1,2,\cdots,n \} \in $ \\ $\mathbb{C}^{m \times m \times p}$, and let $\mathcal{X}_i$ be the difference sequence of $\{\mathcal{Y}_i \}$, i.e., $\mathcal{X}_i \define \mathcal{Y}_i - \mathcal{Y}_{i-1}$ for $i = 1,2,3, \cdots$. If the difference sequence satisfies 
\begin{eqnarray}\label{eq1:cor:7.2}
\mathbb{E}_{i-1} \mathcal{X}_i = 0 \mbox{~~and~~} \mathcal{X}^2_i \preceq  \mathcal{A}_i ~~\mbox{almost surely}, 
\end{eqnarray}
where $i = 1,2,3,\cdots$ and the total varaince $\sigma^2$ is defined as as: $\sigma^2 \define \left\Vert \sum\limits_i^n \mathcal{A}_i^2 \right\Vert$. Then, given a positive real vector $\mathbf{b} \in \mathbb{R}^p$ with $\tilde{j} \define \arg\min\limits_j \{ b_j \}$ and $t \left( \mathcal{Y}_n -  \mathbb{E} \mathcal{Y}_n \right)$ satisfing Eq.~\eqref{eq1:lma: Laplace Transform Method Eigentuple Version} for any $t >0$, we have
\begin{eqnarray}
\mathrm{Pr} \left( \mathbf{d}_{\max}\left( \mathcal{Y}_n -  \mathbb{E} \mathcal{Y}_n \right)\geq \mathbf{b} \right) \leq mp e^{-\frac{b_{\tilde{j}}^2}{8 \sigma^2}}.
\end{eqnarray}
\end{corollary}

Following theorem is the McDiarmid inequality of the maximum eigentuple for the T-product tensor. 

\TensorMcDiarmidEigentuple*

%
\textbf{Proof:}
By the same argument from the proof in Theorem~\ref{thm:Tensor McDiarmid}, this theorem is proved by applying Corollary~\ref{cor:7.2 eigentuple} to the martingale $\{\mathcal{Y}_i\}$.
$\hfill \Box$

\section{Conclusion}\label{sec:Conclusion}.

In Part I paper of this serious work about T-product tensors, we generalize Lapalce transform method and Lieb's concavity theorem from matrices to T-product tensors. In this Part II paper, we apply these techniques to extend the following classical bounds from the scalar to the T-product tensor settings: Chernoff and Bernstein inequalities. The purpose of these probability inequalities tries to identify large-deviation behavior of the extreme eigenvalue and eigentuple of the sums of random T-product tensors. Finally, we also apply these proof techniques developed at this work to study T-product tensor-valued martingales by proving Azuma, Hoeffding, and McDiarmid inequalities under T-product.

\newpage

\bibliographystyle{IEEETran}
\bibliography{Random_TBounds_TProduct_Bib}

\end{document}